\RequirePackage{silence}
\documentclass[12pt]{amsart}
\usepackage{amsmath, amssymb,latexsym,multirow,bm,xcolor}
\usepackage[normalem]{ulem}

\usepackage{mathtools}
\mathtoolsset{showonlyrefs}

\usepackage{aliascnt}

\newtheorem{theo}{{\sc Theorem}}[section]

\newaliascnt{rem}{theo}
\newtheorem{rem}[rem]{{\sc Remark}}
\aliascntresetthe{rem}

\newaliascnt{cor}{theo}

\aliascntresetthe{cor}

\newaliascnt{lem}{theo}
\newtheorem{lem}[lem]{{\sc Lemma}}
\aliascntresetthe{lem}

\newaliascnt{prop}{theo}
\newtheorem{prop}[prop]{{\sc Proposition}}
\aliascntresetthe{prop}

\newaliascnt{defin}{theo}
\newtheorem{defin}[defin]{{\sc Definition}}
\aliascntresetthe{defin}

\newtheorem{maintheo}{{\sc Theorem}}

\newaliascnt{mainprop}{maintheo}
\newtheorem{mainprop}[mainprop]{{\sc Proposition}}
\aliascntresetthe{mainprop}

\newaliascnt{maincor}{maintheo}
\newtheorem{maincor}[maincor]{{\sc Corollary}}
\aliascntresetthe{maincor}

\usepackage[pagebackref]{hyperref}
\hypersetup{
    colorlinks,
    linkcolor={red},
    citecolor={blue},
    urlcolor={blue}
}

\newcommand{\szego}{Szeg\H{o} }
\newcommand{\kahler}{K\"ahler }
\newcommand{\BSj}{Boutet de Monvel-Sj\"ostrand }

\renewcommand{\epsilon}{\varepsilon}
\renewcommand{\phi}{\varphi}

\newcommand{\smallO}{o}

\newcommand{\dashint}{\int}

\newcommand{\bcal}{\mathcal{B}}
\newcommand{\ccal}{\mathcal{C}}
\newcommand{\dcal}{\mathcal{D}}

\newcommand{\fcal}{\mathcal{F}}

\newcommand{\hcal}{\mathcal{H}}

\newcommand{\ocal}{\mathcal{O}}

\newcommand{\ucal}{\mathcal{U}}
\newcommand{\vcal}{\mathcal{V}}

\DeclareMathOperator{\Hess}{Hess}
\DeclareMathOperator{\Ad}{Ad}
\DeclareMathOperator{\Tr}{Tr}

\newcommand{\R}{{\mathbb R}}
\newcommand{\C}{{\mathbb C}}

\newcommand{\Z}{{\mathbb Z}}

\newcommand{\dbar}{\bar\partial}
\newcommand{\ddbar}{\partial\dbar}

\newcommand{\half}{{\frac{1}{2}}}
\newcommand{\Vol}{{\operatorname{Vol}}}

\renewcommand\vec{\bm}

\title[Log-scale equidistribution of zeros]{Log-scale equidistribution
of zeros of quantum ergodic eigensections}

\author{Robert Chang and Steve Zelditch}

\email{hchang@math.northwestern.edu}
\email{zelditch@math.northwestern.edu}

\address{Department of Mathematics, Northwestern University,
Evanston IL,  60208-2730, USA}

\thanks{Research partially supported by NSF grant DMS-1541126}

\begin{document}

\begin{abstract}
Under suitable hypotheses, a symplectic map can be quantized as a sequence of unitary operators acting on the $N$th powers of a positive line bundle over a \kahler manifold. We show that if the symplectic map has polynomial decay of correlations, then there exists a density one subsequence of eigensections whose masses and zeros become equidistributed in balls of logarithmically shrinking radii of lengths $\lvert \log N \rvert^{-\gamma}$ for some constant $\gamma > 0$ independent of $N$.
\end{abstract}

\maketitle

\section{Introduction}

This article is concerned with the equidistribution of masses and of zeros of holomorphic eigensections at the logarithmic scale. Let $(L,h) \rightarrow (M,\omega)$ be a pre-quantum line bundle over a \kahler manifold of complex dimension $m$. In other words, $(L,h)$ is a positive Hermitian line bundle with $c_1(h) = \omega$. Let $(L^N,h^N)$ denote the $N$th tensor power. Under certain quantization conditions (discussed in \autoref{sec:QUANT} and \cite{Z97}), a symplectic map
\begin{equation} \label{CHI}
\chi \colon (M, \omega) \rightarrow (M,\omega), \qquad \chi^* \omega = \omega
\end{equation}
on the base manifold can be quantized as a sequence $\{U_{\chi,N}\}_{N=1}^\infty$ of unitary Fourier integral Toeplitz operators
\begin{equation}
U_{\chi,N} \colon H^0(M,L^N) \to H^0(M, L^N)
\end{equation}
acting on the spaces $H^0(M,L^N)$ of holomorphic sections of $L^N$ with the inner product induced by $h$ (see \autoref{sec:BACKGROUND}).

The eigensections $s_j^N \in H^0(M,L^N)$ of the operators $U_{\chi,N}$ are characterized by
\begin{equation}
U_{\chi,N} s^N_j = e^{i\theta_{N,j}}s_j^N, \qquad 1 \le j \le d_N,
\end{equation}
where $e^{i\theta_{N,j}}$ are eigenphases and $d_N = \dim H^0(M,L^N)$. We write
\begin{equation}
Z_{s^N_j} = \{z \in M \colon s^N_j(z) = 0\} \quad \text{and} \quad \big[Z_{s^N_j}\big] = \frac{\sqrt{-1}}{\pi} \ddbar \log \|s^N_j(z)\|^2_{h^N} + N\omega
\end{equation}
for the zero set of $s^N_j$ and the current of integration over the zero set of $s^N_j$, respectively (cf.\ \eqref{PL}). With the ergodicity assumption on $\chi$, \cite{Z97} proved  that the eigensections of the quantum maps $U_{\chi,N}$ are quantum ergodic. Moreover, \cite{NV98} and \cite{ShZ99} (see also \cite{R05} for the modular surface setting) proved that the zeros of `almost all' quantum ergodic eigensections are asymptotically equidistributed with respect to the \kahler volume form: There exists a subsequence $\Gamma
\subset \{(N, j) \colon \text{$N \geq 1$, $j = 1, \dots, d_N$}\}$ of density one for which
\begin{equation}\label{eqn:EQUI ZERO}
\lim_{\substack{(N,j) \in \Gamma \\ N\rightarrow\infty}} \int_M  f(z) \left[\frac{1}{N} Z_{s_j^N} \right] \wedge \omega^{m-1} = \int_M f \frac{\omega^m}{m!}, \qquad f \in C^\infty(M).
\end{equation}

\subsection{Statement of main results}\label{STATEMENT} 

Recall that $m = \dim_\C M$. We fix a logarithmic scale $\epsilon_N$ depending on parameter $\gamma$:
\begin{equation}\label{eqn:EPSILON}
\epsilon_N := \lvert\log N \rvert^{-\gamma} \quad \text{for some constant $0 < \gamma < \frac{1}{6m}$ independent of $N$.}
\end{equation}
The main purpose of this paper is to  show  (with additional assumptions on $\chi$, described below) that the equidistribution result \eqref{eqn:EQUI ZERO} holds with the domain of integration $M$ replaced by any ball $B(p, \epsilon_N)$ centered at $p \in M$ with radius $\epsilon_N = \lvert \log N \rvert^{-\gamma}$ for any $\gamma < (6m)^{-1}$. This is what is meant by ``equidistribution of zeros at the logarithmic scale.''

To obtain this log-scale improvement, we use two dynamical properties of $\chi$:
\begin{itemize}
\item For $T \in \mathbb{Z}$, let $\chi^T$ denote the $T$-fold iterate of $\chi$ (or of its inverse $\chi^{-1}$, depending on the sign of $T$). By the chain-rule $\chi$ satisfies the exponential growth estimate
\begin{equation}\label{eqn:EXPGROWTH}
\|\chi^T\|_{C^2} = \ocal(e^{ |T|\delta_0}) \quad \text{for some fixed constant $\delta_0 > 0$ independent of $T$.} 
\end{equation}
In particular, if $\chi$ lifts to a contact transformation $\tilde{\chi}$ on the unit co-disk bundle $X \rightarrow M$ (see \autoref{sec:QUANT}), then  $\|F \circ \tilde{\chi}^\ell\|_{C^2}^2 = \ocal_F(e^{2|T| \delta_0})$ for any $F \in C^{\infty}(X)$.

\item We  assume that $\chi$ has sufficiently fast decay of correlations. Namely, that there exist constants $0 < \beta < 1$,  $c_1 > 0$, and $c_2 = c_2(\beta)> 1$ such that\footnote{Even though an exponential decay rate (i.e., with $(1 +\lvert T \rvert)^{-c_2}$ replaced by $e^{ -c_2 \lvert T \rvert}$) is often assumed in the literature, much less is necessary for the proof; this was also noted in \cite{Sc06}.}
\begin{equation}\label{eqn:DECAYOFCORR}
\left\lvert \int_M (g \circ \chi^T)f  \,dV  - \int_M f\,dV \int_M g\,dV\right\rvert \le c_1  (1 +\lvert T \rvert)^{-c_2} \|f\|_{C^{0,\beta}}\|g\|_{C^{0,\beta}}
\end{equation}
for all $f,g \in C^{0,\beta}(M)$. Thus,  $\chi$ is mixing and hence ergodic. Here, $dV$ is the normalized volume form \eqref{NORMVOL}.
\end{itemize}
The explicit error estimate in Egorov's theorem for Toeplitz operators (\autoref{prop:EGOROV}, proved in \autoref{sec:EGOROV}) relies on assumption \eqref{eqn:EXPGROWTH}. Assumption \eqref{eqn:DECAYOFCORR} is used in the proof of logarithmic decay of quantum variances (\autoref{theo:LOGQE}) in \autoref{sec:LOGQE}.

\subsubsection{Log-scale equidistribution of zeros}

The log-scale equidistribution of zeros states that zeros in balls of radii $\epsilon_N$ are uniformly distributed with respect to the volume form \eqref{NORMVOL}. It
is simplest to state the result by dilating such shrinking balls by $\epsilon_N^{-1}$ back to
a fixed reference ball of radius $1$. 
In a local  \kahler normal  coordinate chart $(U, z)$ with
$z = 0$ at $p$, define local dilation maps 
\begin{equation}\label{Dep}
D^p_{\epsilon} \colon B(p, 1) \to B(p, \epsilon), \qquad D_{\epsilon} z = \epsilon z.
\end{equation}
 Here we abuse notation by writing $B(p,1)$ when we mean the image of the metric unit ball centered at $p$ in the local coordinate chart based at $p$.
  The inverse dilation is defined by 
\begin{equation}
(D^p_{\epsilon})^{-1} \colon B(p, \epsilon) \to B(p,1).
\end{equation}

\begin{rem} 
We recall that  \kahler normal coordinates $z_1, \cdots, z_m$ centered at point $z_0$ are holomorphic coordinates in which $z_0$ has coordinates $0 \in \C^m$, and 
\[ \omega(z) = i \sum_{j=1}^m d z_j \wedge d \bar z_j + O(|z|^2). \]
We may also choose a local reference frame $e_L$ of the line bundle in a neighborhood of $z_0$, such that the induced \kahler potential $\varphi$ takes the form
\[ \varphi(z) = |z|^2 + O(|z|^3). \]
We refer to \cite{GriffithsHarris78} for background.
\end{rem}

Let $D_{\epsilon}^{p*}$ be the corresponding pullback operator on forms. 
For simplicity of notation we denote the pullback $(D_{\epsilon}^p)^{* -1}$ of the inverse dilation by
$D_{\epsilon *}^p$ so that
\begin{equation}
D_{\epsilon *}^p \colon \dcal^{m-1,m-1}(B(p,1)) \to \dcal^{m-1,m-1}(B(p, \epsilon)),
\end{equation}
where $\dcal^{m-1,m-1}$ denotes the space of
 compactly supported smooth $(m-1, m-1)$ test forms. In particular, for  $\eta \in  \dcal^{m-1,m-1}(B(p, 1))$, we  have 
\begin{equation}
\int_{B(p, \epsilon)}  D_{\epsilon *}^p \eta \wedge \frac{1}{N} \big[Z_{s_j^N}\big] = \int_{B(p,1)}  \left(\eta \wedge \frac{1}{N} D_{\epsilon}^{p*} \big[Z_{s_j^N}\big] \right).
\end{equation}

\begin{maintheo}\label{theo:LOGZERO}
Let $(L, h) \rightarrow (M, \omega)$ be a pre-quantum line bundle. Let $\chi$ satisfy \eqref{eqn:EXPGROWTH} and \eqref{eqn:DECAYOFCORR}.
Let  $\{s_1^N, \dots, s_{d_N}^N\}$ be an
orthonormal basis of eigensections of $U_{\chi, N} \colon H^0(M, L^{N}) \rightarrow  H^0(M, L^{N})$. Then, for every $0 < \gamma < (6m)^{-1}$ and $\epsilon_N = \lvert \log N \rvert^{-\gamma}$, there
exists a full density subsequence $\Gamma  \subset \{(N, j) \colon j = 1, \dots, d_N\}$ such that for every $p \in M$,
\begin{equation}
\frac{1}{N\epsilon_N^2} D_{\epsilon_N}^{p*} \big[Z_{s_j^N}\big]  \xrightharpoonup{\Gamma \ni (N,j) \rightarrow \infty} \omega^p_0 \quad \text{in the weak sense of currents on $B(p,1)$,}
\end{equation}
where $\omega^p_0 = \frac{\sqrt{-1}}{2\pi} \ddbar \log |z|^2$ is the flat \kahler form  in \kahler normal coordinates based
at $p$.
\end{maintheo}

\begin{rem}
The weak convergence statement in \autoref{theo:LOGZERO} means that for every test form $\eta \in \dcal^{m-1,m-1}(B(p,1))$, one has
\begin{equation}
\int_{B(p,1)} \left(\eta \wedge\frac{1}{ N\epsilon_N^2} D_{\epsilon_N}^{p*} \big[Z_{s_j^N}\big]  \right) = \int_{B(p,1)} \eta \wedge
\omega_0^p  + \smallO(1).
\end{equation}
\end{rem}
The key ingredients of the proof are the log-scale mass comparison result (\autoref{theo:LOGMASS}), the Poincar\'e-Lelong formula \eqref{PL} and
compactness results on logarithms of scaled sections. This equidistribution result should be compared to 
Lester-Matom\"aki-{Radziwi\l\l}  \cite[Theorem~1.1]{LMR15} for a sequence $\{f_k\}$ of Hecke modular cusp forms of weight $k$. 
They proved 
that for a certain $\delta > 0$,
\begin{equation}
\frac{\# \{ z \in B(z_0, r) \colon  f_k(z) = 0\}}{\# Z_{f_k}} = \frac{3}{\pi} \int_{B(z_0, r)} \frac{dx dy}{y^2}
 + \ocal \left( r (\log k)^{- \delta + \epsilon} \right)
\end{equation}
when $r \geq (\log k)^{-\delta/2 + \epsilon}$.
This  is a quantum \emph{unique} ergodicity result in that it  holds for the entire orthonormal basis of Hecke eigenforms, whereas we discard a density zero subsequence of eigensections because we work in the more general setting of dynamical Toeplitz quantizations of  quantizable ergodic symplectic maps.

\subsubsection{Log-scale equidistribution of mass}

The equidistribution result of \autoref{theo:LOGZERO} is based on 
log-scale volume comparison theorems similar to  those of \cite[Lemma~3.1]{HR16} and  \cite[Corollary~1.9]{H15}.

\begin{maintheo}[Log-scale equidistribution of masses]\label{theo:LOGMASS}
Assume the hypotheses of \autoref{theo:LOGZERO}. Then, given any $0 < \gamma' < (6m)^{-1}$ and $\epsilon_N' = \lvert \log N \rvert^{-\gamma'}$ as defined by \eqref{eqn:EPSILON}, there exist a full density subsequence $\Gamma$  and constants $C_1, C_2$  uniform in $p \in M$ and independent of $N$ such that
\begin{equation}
 C_1 \frac{\Vol(B(p, \epsilon_N'))}{\Vol(M)} \le \int_{B(p, \epsilon_N')} \| s_j^N\|_{h^N}^2 \, dV \le C_2 \frac{\Vol(B(p, \epsilon_N'))}{\Vol(M)} \quad \text{as $\Gamma \ni (N,j) \rightarrow \infty$}.
\end{equation}
Here, $dV$ is the normalized volume form \eqref{NORMVOL}.
\end{maintheo}
There is no need to put primes on $\gamma$ or $\epsilon_N$ in the statement above, but we do so to foreshadow that in the proof of \autoref{theo:LOGZERO}, the result of \autoref{theo:LOGMASS} is applied with $\gamma < \gamma'$ and $\epsilon_N' < \epsilon_N$. The comparison (as opposed to asymptotic) result on log-scale mass equidistribution is sufficient for deriving equidistribution of zeros at a slightly larger logarithmic scale. In fact, only the lower bound is used, and the bound itself is much stronger than necessary for the proof. 

\autoref{theo:LOGMASS} is based on a quantitative quantum variance estimate (\autoref{theo:LOGQE}) in the holomorphic setting. Before stating the estimate, we record here another one of its corollaries, which is analogous to \cite[Corollary~1.8]{H15}.

\begin{mainprop}\label{cor:CESARODENSITYONE}
Assume the hypotheses of \autoref{theo:LOGZERO}. Fix $z_0 \in M$. Then, given any $0 < \gamma < (4m)^{-1}$ and $\epsilon_N$ as defined by \eqref{eqn:EPSILON}, there exists a subsequence $\Gamma_{z_0} \subset \{(N,j)\}$ of density one such that
\begin{equation}\label{eqn:QEFIXED}
\int_{B(z_0,\epsilon_N)} \|s_j^N\|_{h^N}^2\,dV = \frac{\Vol(B(z_0,\epsilon_N))}{\Vol(M)} + \smallO(\lvert \log N \rvert^{-2m\gamma}).
\end{equation}
Here, $dV$ is the normalized volume form \eqref{NORMVOL}.
\end{mainprop}
Recall  $\dim_\C M = m$, so $\frac{\Vol(B(z_0,\epsilon_N))}{\Vol(M)} = C(M,g) \epsilon_N^{2m} = C(M,g) \lvert \log N \rvert^{-2m\gamma}$. The differences between \autoref{cor:CESARODENSITYONE} and \autoref{theo:LOGMASS} are that the former is an asymptotic result for a fixed base point, whereas the latter is a comparison result that holds for all points in $M$. Moreover, in the former case the range of values that $\gamma$ can take is improved. \autoref{cor:CESARODENSITYONE} is not used in proving  \autoref{theo:LOGZERO} or \autoref{theo:LOGMASS}.

\subsubsection{Log-scale quantum ergodicity}\label{sec:QEintro}
By the quantum variance associated to $f$ we mean the quantity
\begin{equation}\label{eqn:VARIANCE}
\vcal_N(f) := \frac{1}{d_N} \sum_{j=1}^{d_N} \left\lvert \,\int_M f(z)\|s^N_j\|_{h^N}^2\,dV - \dashint_M  f\,dV\right\rvert^2 \quad \text{for $f \in C^{\infty}(M)$.}
\end{equation}
Here, $dV$ is the normalized volume form (7). Thanks to Egorov's theorem for Toeplitz operators (\autoref{prop:EGOROV}, proved in \autoref{sec:EGOROV}) and the decay of correlations assumption \eqref{eqn:DECAYOFCORR}, we show the quantum variance has a logarithmic decay rate when $f \in C^{\infty}(M)$:
\begin{maintheo}[Logarithmic decay of quantum variances]\label{theo:LOGQE}
Assume the hypotheses of \autoref{theo:LOGZERO}. Then, there exists a constant $\kappa_0 > 0$ independent of $N$  such that for every $0 < \beta < 1$ and for every $f \in C^{2}(M)$,
\begin{equation}\label{C2}
\vcal_N(f) = \ocal\left(\frac{\|f\|_{C^{0,\beta}}^2}{\log N}\right) +\ocal\left(\frac{\|f\|_{C^2}^2\lvert\log N\rvert^2}{N^{\frac{1}{2}}}\right)
+ \ocal\left(\frac{\|f\|_{C^{0,\beta}}^2}{N\log N}\right),
\end{equation}
where $\| \cdot \|_{C^{0,\beta}}$ is the $\beta$-H\"older norm.
\end{maintheo}

We specialize to the following logarithmically dilated symbols. In \kahler normal coordinates, let $f_{z_0} \in C_0^{\infty} (B(z_0,2), \R)$  be a smooth cut-off function  that is equal to 1 on $B(z_0,1)$, vanishes outside of $B(z_0,2)$ and satisfies $0 \le f_{z_0} \leq 1$. For ``small-scale quantum ergodicity,'' we work with locally dilated  symbols (cf.\ \eqref{Dep}) of the form 
\begin{equation}\label{eqn:LOGSCALEQE}
f_{z_0,\epsilon}(z) := D^{z_0}_{\epsilon *}f_{z_0} (z)=  f\Big( \frac{z}{\epsilon}\Big) \in C_0^{\infty} (B(z_0,2 \epsilon), \R), \quad \text{where $z_0 \in M$ and $\epsilon > 0$}.
\end{equation}
Then set $\epsilon = \epsilon_N$. It follows from \autoref{theo:LOGQE} that, to leading order in $N$, the quantum variance associated to such symbols have the estimate
\begin{align}
\vcal_N(f_{z_0,\epsilon_N}) &= \ocal(\|f_{z_0}\|_{C^{0,\beta}}^2\lvert \log N \rvert^{2\gamma\beta - 1}).
\end{align}
Since $0 < \beta < 1$ and $\gamma < \frac{1}{6m}$, we have $2\gamma\beta - 1 < 0$. Since the second term is smaller than the first,
we obtain:
\begin{maincor}[Log-scale quantum variance estimates]\label{cor:LOGQE}
Let $\epsilon_N$ be as defined in \eqref{eqn:EPSILON}. Under the same hypotheses as in \autoref{theo:LOGZERO}, we have
\begin{equation}\label{eqn:SMALLSCALEQE}
\vcal_N(f_{z_0,\epsilon_N}) = \ocal(\|f_{z_0}\|_{C^{0,\beta}}^2\lvert \log N \rvert^{2\gamma\beta - 1}),
\end{equation}
where the error estimate is uniform in $z_0$.

\end{maincor}

Following the arguments of \cite{HR16} and \cite{H15}, an application of \autoref{cor:LOGQE} and a covering argument together imply
\autoref{theo:LOGMASS}.

\subsection{Further results}
The results of this paper are the \kahler analogue of the small-scale quantum ergodicity results in the Riemannian setting proved in \cite{HR16,H15}. Specializing to the torus $\mathbb{T}^d = \R^d/2\pi \Z^d$, \cite[Theorem~1.1]{LR} proves the stronger uniform mass distribution result
\begin{equation}
\lim_{n \rightarrow \infty} \sup_{B(y,r) \in \bcal_n}\left\lvert \frac{1}{\Vol(B(y,r))} \int_{B(y,r)} \lvert \psi_n \rvert^2\,dx - 1 \,\right\rvert = 0,
\end{equation}
for a density one subsequence of eigenfunctions $(\Delta - \lambda_j)\psi_j = 0$. Here, the supremum is taken over the set $\bcal_n$ of balls $B(y,r) \subset \mathbb{T}^d$ of radii $r > \lambda_n^{-1/(2d-2) + o(1)}$.

For Hecke modular eigenforms, it is proved in  \cite[Theorem~1.5]{LMR15} that (in the notation defined above)
\begin{equation}
\sup_{{\mathcal R} \subset \fcal}  \left| \int_{{\mathcal R}} y^k |f_k(z)|^2 \frac{dx dy}{y^2} - \frac{3}{\pi} \int_{{\mathcal R}}
\frac{dx dy}{y^2} \right| \leq C_{\epsilon}  (\log k)^{- \delta + \epsilon},
\end{equation}
where the supremum is taken over all rectangles $\mathcal{R}$ with sides parallel to the $x$- and $y$-axis. This is a stronger result because it is valid for
all Hecke eigenforms and because  the supremum is taken over rectangles of any size rather than over rectangles of `radii'
$\epsilon_k = \lvert\log k\rvert^{-\gamma}$.

In the \kahler setting, \cite{ShZ99} proves equidistribution of zeros (not at the logarithmic scale) for random orthonormal
bases of $H^0(M, L^N)$ as well as for eigensections of quantized ergodic
symplectic maps. It is probable that \autoref{theo:LOGZERO} can also
be generalized to random orthonormal bases. This is work in progress of
the first author.

\subsection{Existence of quantizable ergodic symplectic diffeomorphisms}

An obvious question is whether quantizable ergodic symplectic diffeomorphisms satisfying the decay of correlations condition \eqref{eqn:DECAYOFCORR} 
exist on a given  \kahler manifold. (Any diffeomorphism satisfies the exponential growth estimate \eqref{eqn:EXPGROWTH} automatically.)
There seem to exist  few studies of ergodic symplectic dynamics in 
dimensions $> 2$.   After consulting with  several experts in the field, we  give a brief summary of the examples that we are aware of. 

The simplest and most-studied examples are
 hyperbolic symplectic toral automorphisms induced by an element
of $Sp(2n, \Z)$   and small perturbations
 of such automorphisms  (see \cite{Z97,Kel10} for their Toeplitz  quantizations). More generally, any  hyperbolic or Anosov symplectic diffeomorphism satisfies the assumptions.
 There  is a quantization condition, but  as explained in \cite{FT15}, it is always satisfied  if one tensors with a flat line bundle and
modifies the contact form.

Most studies of smooth ergodic maps concern  volume  preserving 
diffeomorphisms. Studies of ergodic symplectic diffeomorphisms on manifolds other than tori are rare except in the dimension two. In that dimension, ergodic (indeed, Bernoulli) symplectic diffeomorphisms of surfaces
of any genus exist (see~\cite{Ka79} and Theorem~1.26 of  \cite{BP13}). As mentioned above,
they are quantizable.  We also mention  that    pseudo-Anosov diffeomorphisms are singular ergodic symplectic 
diffeomorphisms which are smooth  away from a finite number of singular points. They act hyperbolically with respect to  two transverse (singular) measured foliations. 
Since they are singular, our techniques do not apply directly but it is plausible that they can be modified by suitably cutting off singular points. These examples may turn out to be the most explicitly computable ones on surfaces other than tori and are very likely to satisfy all the conditions of this article.

In higher dimensions, Anosov diffeomorphisms have been studied on 
certain types of nilmanifolds in addition to tori (see \cite{DeV11}). Partially hyperbolic symplectic diffeomorphisms are studied in \cite{M16}.
There are further  partially hyperbolic examples obtained by perturbation.
 As explained to the authors by A. Wilkinson, a symplectic toral automorphism (or any partially
hyperbolic symplectic diffeomorphism)   can be perturbed to produce 
a symplectic diffeomorphism which is stably accessible (see \cite{DW03}). Moreover, if the original map is ``center bunched,'' then the perturbed map 
is stably ergodic (see \cite{BW10}).   These examples are additional to the
usual Anosov diffeomorphisms of tori and their perturbations. We refer to these articles for the definitions
and further discussion. 

\subsection{Acknowledgments} We thank H. Hezari for pointing some errors and gaps in the earlier version, and for suggesting corrections. We also thank  F. Faure,  G. Riviere and A. Wilkinson for useful comments and references on the dynamical aspects. Finally, we thank the referees for their detailed and helpful comments that led to significantly improvements of the paper.

\section{Background}\label{sec:BACKGROUND}

\subsection{Complex geometry}\label{CG} 

We follow the notation used in \cite{ShZ99,Z97,Z98} and refer there for further discussion. Let $(M,\omega)$ be a compact \kahler manifold of dimension $\dim_\C M = m$. Let $(L,h) \rightarrow (M,\omega)$ be a pre-quantum line bundle. In other words, $L$ is a ample Hermitian line bundle endowed with a smooth metric $h$ whose curvature form $c_1(h)$ is strictly positive with $c_1(h) = \omega$. If $e_L$ is a nonvanishing local holomorphic frame for $L$ over an open set $U \subset M$, then
\begin{equation}
c_1(h)=- \frac{\sqrt{-1}}{\pi}\ddbar \log \|e_L\|_h,
\end{equation} 
where $ \|e_L\|_h := h(e_L, e_L)^{1/2}$ denotes the $h$-norm of $e_L$.

The curvature form $c_1(h)$ is a representative of the first Chern class $c_1(L) \in H^2(M, \R)$ of the line bundle $L$. Since $c_1(h) = \omega$, we have $\int_M \omega^m = c_1(L)^m \in \Z^+$. We normalize the volume form by this quantity so that $M$ has unit volume:
\begin{equation}\label{NORMVOL}
dV := \frac{\omega^m}{c_1(L)^m}.
\end{equation}

We work with spaces $H^0(M,L^N)$ of holomorphic sections $s^N$  of $L^N$. (The superscript on $s$ indexes
the degree and does not mean the $N$th power.) These are finite dimensional Hilbert spaces of dimensions
\begin{equation}
d_N := \dim H^0(M,L^N) \sim \frac{c_1(L)^m}{m!}N^m \quad \text{as $N \rightarrow \infty$.}
\end{equation}
The Hermitian metric $h^N$ and the inner product structure on $H^0(M,L^N)$ are tensor powers of  the metric $h$ on $L$:
\begin{equation}
\begin{dcases}
 \|s^{\otimes N}(z)\|_{h^N} := \|s(z)\|_h^N & s \in H^0(M,L),\\
\langle s_1^N,s_2^N \rangle = \int_M h^N(s_1^N(z),s_2^N(z))\,dV & s_1^N, s_2^N \in H^0(M,L^N).
\end{dcases}
\end{equation}

Given a holomorphic section $s^N \in H^0(M,L^N)$, we denote by $\big[Z_{s^N}\big]$ its current of integration over the zero divisor of $s^N$. In a local frame $e_L^N$ for $L^N$, we can write $s^N=f^{(N)} e_L^N$ with
$f^{(N)}$ a holomorphic function. Let $g(z) := \|e_L(z)\|_{h}^2 = e^{-\phi(z)}$ where $\phi$ is the \kahler potential, then $\|e_L^{N}(z)\|_{h^N}^2 = g(z)^N$ and $\|s^N\|_{h^N}^2 = |f^{(N)}|^2 g^N$. The Poincar\'e-Lelong formula states that
\begin{equation} \label{PL}
\big[Z_{s^N}\big] =\frac{\sqrt{-1}}{\pi} \ddbar \log|f^{(N)}| =
\frac{\sqrt{-1}}{\pi} \ddbar \log \|s^N\|_{h^N} + N\omega.
\end{equation}

\subsection{Hardy space of CR holomorphic functions}

Let $(L^*, h^*)$ be the dual line bundle to $L \rightarrow M$. Thanks to the positivity of $c_1(h)$, the unit co-disk bundle $D^* \subset L^*$ relative to the dual metric $h^*$ is a strictly pseudoconvex domain whose boundary
\begin{equation}
X := \partial D^* = \{v \in L^* \colon h^*(v,v) = 1\} \subset L^*
\end{equation}
is a CR manifold. The Hardy space $H^2(X)$ is the space of square integrable CR functions on $X$, or equivalently the space of boundary values of holomorphic functions on the unit disk bundle with finite $L^2(X)$ norm.

We introduce a defining function $\rho$ for $X$, which will be featured in the \BSj parametrix. We write points in the co-disk bundle as $x = (z,\lambda e_L^*(z))$, where $\lambda \le 1$ and $e_L^*(z)$ is a normalized dual frame centered at $z \in M$. Define
\begin{equation}\label{eqn:DEFININGFCTN}
\rho \colon D^* \rightarrow \mathbb{R}, \quad \text{$\rho(z, \lambda e_L^*(z)) = 1 - \lvert \lambda \rvert^2 e^{-\phi(z)}$ where $\phi$ is the \kahler potential}.
\end{equation}
Then $\rho$ is a defining function for $X$ satisfying
\begin{itemize}
\item $\rho$ is defined in a neighborhood of $X$;

\item $\rho > 0$ in $D^*$;

\item $\rho = 0$ on $X$;

\item $d\rho \neq 0$ near $X$.
\end{itemize}
We define the contact form $$\alpha = d^c \rho |_X. $$

Let $r_\theta$ be the natural circle action on $X$, that is, $r_\theta x = e^{i\theta} x$ for $x \in X$. Note that a section $s \in H^0(M,L)$ determines an equivariant function $\hat{s}$ on $L^*$ by the rule
\begin{equation}
\hat{s}(z,\lambda) = (\lambda, s(z)), \qquad \text{$z \in M$, $\lambda \in L^*_z$}.
\end{equation}
It is easy to verify restricting $\hat{s}$ to $X$ yields $\hat{s}(r_\theta x) = e^{i\theta}\hat{s}(x)$. Conversely, a section $s^N \in H^0(M,L^N)$ determines an equivariant function $\hat{s}^N$ on $L^*$ whose restriction to $X$ satisfies $\hat{s}{(N)}(r_\theta x) = e^{iN\theta}\hat{s}^N(x)$. The map $s^N \mapsto \hat{s}^N$ is in fact a unitary equivalence between the space $H^0(M,L^N)$ of holomorphic sections and the weight spaces
\begin{equation*}
H^2_N(X) := \left\{F \in H^2(X) \colon F(r_\theta x) = e^{iN\theta}F(x)\right\} \quad \text{with} \quad H^2(X) = \bigoplus_{N \ge 0} H^2_N(X).
\end{equation*} 

The   \szego projector is the orthogonal projection $$\Pi \colon L^2(X) \rightarrow H^2(X)$$
and its Fourier components are denoted by
\begin{equation}
\Pi_{h^N} \colon L^2(X) \rightarrow H_N^2(X).
\end{equation}

\subsection{Quantization of symplectic maps}\label{sec:QUANT} 

We use the dynamical Toeplitz quantization method of \cite{Z97}. A symplectic map $\chi \colon M \rightarrow M$ is quantizable if and only if it lifts to a connection-preserving contact transformation $\tilde{\chi} \colon X \rightarrow X$, that is, $\tilde{\chi}^* \alpha = \alpha$. Denote by
\begin{equation}
T_{\tilde{\chi}} \colon L^2(X) \rightarrow L^2(X), \qquad T_{\tilde{\chi}} F = F \circ \tilde{\chi}
\end{equation}
the pre-composition by the lift $\tilde{\chi}$.  Note that $\tilde{\chi}$ commutes with the natural circle action $r_\theta$ on $X$, and $\|\tilde{\chi}\|_{C^2(X)} = c \cdot \|\chi\|_{C^2(M)}$ for some constant $c$.

The quantization of a quantizable map $\chi$ is defined to be a unitary Fourier integral operator
\begin{equation}\label{eqn:QUANTIZE2}
U_\chi := \Pi  \sigma T_{\tilde{\chi}}\Pi \colon H^2(X) \rightarrow H^2(X).
\end{equation}
Here, $\sigma$ is a zeroth order symbol that makes the operator $U_\chi$ defined by \eqref{eqn:QUANTIZE1} unitary. Its  existence is guaranteed by the construction in \cite{Z97}).   We emphasize again that $T_\chi$ denotes translation by the lifted map; such translation is not well-defined on the base because it does not preserve the line bundle.

Under the identification $H^2(X) = \bigoplus_{N \ge 0} H^2_N(X)$, $U_\chi$ decomposes into a sequence of unitary Fourier integral operators $U_{\chi,N}$ defined by
\begin{equation}\label{eqn:QUANTIZE1}
U_{\chi,N} := \Pi_{h^N}  \sigma_N T_{\tilde{\chi}} \Pi_{h^N} \colon  H^2_N(X) \rightarrow  H^2_N(X).
\end{equation}
Here, $\sigma_N$ is a zeroth order symbol making $U_{\chi,N}$ unitary. The Fourier coefficients $\Pi_{h^N}$ have an explicit parametrix given in \eqref{eqn:PARAMETRIX}.

\subsection{\BSj parametrix for the \texorpdfstring{\szego}{Szego} projector}\label{sec:BSj}

In preparation for the proof of Egorov's theorem for Toeplitz operators (\autoref{prop:EGOROV}), we briefly recall the \BSj parametrix for the \szego kernel. Let $\Pi(x,y)$ denote the kernel of the \szego projector $\Pi$ in \eqref{eqn:QUANTIZE2}, that is,
\begin{equation}
\Pi F(x) = \int_X \Pi(x,y)F(y)\,dV(y) \quad \text{for all $F \in L^2(X)$}.
\end{equation}
It is proved in \cite{BS75} that $\Pi$ is a complex Fourier integral operator of positive type. Near the diagonal, there is a parametrix of the form
\begin{equation}
\Pi(x,y) \sim \int_0^\infty e^{it \psi(x,y)} s(x,y,t)\,dt,
\end{equation}
where
\begin{equation}
s(x,y,t) \sim \sum_{n =0}^\infty t^{m - n} s_n(x,y)
\end{equation}
belongs to the symbol class $S^m(X \times X \times \R_{\ge 0})$ and $\psi \in C^\infty(D^* \times D^*)$ is a complex phase of positive type. (Recall that $D^*$ stands for the unit co-disk bundle, of which $X$ is the boundary.)

The phase function $\psi$ is obtained as the almost-analytic continuation of the defining function $\rho$ in \eqref{eqn:DEFININGFCTN}. Explicitly, for $x_j = (z_j,\lambda_j e_L^*(z_j)) \in D^*$, we have
\begin{equation*}
\psi(x_1,x_2) = \frac{1}{i} \bigg(1 - \lambda_1\bar{\lambda}_2 e^{-\frac{\phi(z_1)}{2} - \frac{\phi(z_2)}{2} + \phi(z_1,\bar{z}_2)}\bigg),
\end{equation*}
where $\phi(z_1,\bar{z}_2)$ is obtained from the \kahler potential $\phi$ by writing $\phi(z_1) = \phi(z_1,\bar{z}_1)$ on the diagonal of $M \times \overline{M}$ and extending to a neighborhood of the diagonal. When the metric is real analytic the extension is analytic; in the general $C^{\infty}$ case it is almost-analytic.  If we assume in addition that $x_j \in X$ lie on the co-circle bundle, then $\lambda_j = e^{i\tau_j}$ is uni-modular, whence $x_j = (z_j, \tau_j)$ and
\begin{equation}\label{eqn:SMALL PHASE}
\psi(x_1,x_2) = \psi(z_1,\tau_1,z_2,\tau_2) = \frac{1}{i} \bigg(1 - e^{-\frac{\phi(z_1)}{2} - \frac{\phi(z_2)}{2} + \phi(z_1,\bar{z}_2)}e^{i(\tau_1 - \tau_2)}\bigg) \quad \text{on $X \times X$}.
\end{equation}

The kernels of the partial \szego projectors $\Pi_{h^N}$ in \eqref{eqn:QUANTIZE1} are the Fourier coefficients of $\Pi(x,y)$:
\begin{align}\label{eqn:PARAMETRIX}
\Pi_{h^N}(x,y) &= \int_0^\infty \! \int_{S^1} e^{-iN\theta} e^{it\psi(r_\theta x, y)} s(r_\theta x,y,t)\,d\theta dt\\ \notag
&=N\int_0^\infty\!\int_{S^1} e^{iN[-\theta + t \psi(r_\theta x, y)]} s(r_\theta x,y,Nt)\,d\theta dt,
\end{align}
where the second line follows from a change of variable $t \mapsto Nt$.

\subsection{Off-diagonal estimates and scaling asymptotics}\label{ODSECT} 

We will be using two off-diagonal estimates for the lifted \szego kernel on
$X \times X$. Again, write $x_j = (z_j,\tau_j)$ for points in the co-circle bundle $X$. Let $d(z,w)$ be the distance with respect to the \kahler metric on $M$. 

The first is an Agmon-type estimate giving global off-diagonal bounds:
\begin{equation}\label{offdiag}
\lvert \Pi_{h^N}(x_1, x_2)| \le A_1N^m e^{-A_2\sqrt{N}d(z_1,z_2)} \quad \text{for constants $A_1,A_2$ independent of $N,x_1,x_2$}
\end{equation}
due to Lindholm \cite{Lindholm01}, Delin \cite{Delin98} and others. The second is a near diagonal Gaussian decay estimate: There exists $A_3 < 1$ independent of $N,x_1,x_2$ such that 
\begin{equation}\label{neardiag}
\lvert\Pi_{h^N}(x_1,x_2)\rvert \le \bigg(\frac 1{\pi^m}
+o(1)\bigg) {N^m}e^{-\frac {1-A_3} 2 N
d(z_1,z_2)^2}+O(N^{-\infty}) \quad \text{whenever  $d(z, w) \leq N^{-\frac{1}{3}}$}. \end{equation} 
We refer to \cite{ShZ99,ShZ02,MaMa12} for background and references.

We further use near off-diagonal scaling asymptotics from \cite{ShZ02, LuSh15}.
At
each $z \in M$ there is an osculating Bargmann-Fock or Heisenberg model associated to $(T_z M, J_z, h_z)$.
Let $(u, \theta_1, v, \theta_2)$ be linear coordinates on $T_z M \times S^1 \times T_z M
\times S^1$. The model Heisenberg \szego kernel on the tangent space is denoted by 
\begin{equation} \label{TANGENT} 
\Pi^{T_z M}_{h_z, J_z}(u, \theta_1,  v, \theta_2) : L^2(T_z M) \to \hcal(T_z M, J_z, h_z)
= \hcal_J.
\end{equation}
We recall that the semi-classical \szego kernels of the Heisenberg group have the form
\begin{equation}\label{szegoheisenbergintro}
\Pi_{h^N}^\mathbf{H}(x_1,x_2)  =\frac{1}{\pi^m} N^m
e^{i N (\tau_1 - \tau_2 )} e^{ N(z_1 \cdot \bar{z}_2 -\frac{1}{2} \lvert z_1 \rvert^2
-\frac{1}{2}\lvert z_2 \rvert^2) }.\end{equation}

In \cite{LuSh15} the notion of K-coordinates is introduced, refining the notion of Heisenberg coordinates in \cite{ShZ02}. These are K\"ahler-type coordinates in which \eqref{TANGENT} equals \eqref{szegoheisenbergintro}
to leading order (up to rescaling):
\begin{equation} \label{SCALINGB} \Pi^{T_z M}_{h_z, J_z}(u, \theta_1,  v, \theta_2)  = \pi^{-m} e^{i (\theta_1 - \theta_2)} e^{u\cdot \bar{v} - \half (|u|^2 + |v|^2)}=  \pi^{-m} e^{i (\theta_1 - \theta_2)} e^{i \Im (u\cdot \bar{v}) - \frac{1}{2} |u -v|^2}
\end{equation}

The lifted \szego kernel is shown in \cite{ShZ02}   and in Theorem~2.3 of \cite{LuSh15} to have the following scaling asymptotics.
\begin{theo} \label{SHLU} Fix $P_0\in M$ and choose a K-frame centered at $P_0$. Then, identifying coordinates $(z_1,\tau_1,z_2,\tau_2)$ on $X^2$ with coordinates $(u,\theta_1,v,\theta_2)$ on $(T_zM \times S^1)^2$, we have
\begin{multline} \label{SCALING} N^{-m} \Pi_{h^N}\bigg(\frac{u}{\sqrt{N}}, \frac{\theta_1}{N}, \frac{v}{\sqrt{N}}, \frac{\theta_2}{N}\bigg)\\
 = \Pi^{T_z M}_{h_z, J_z}(u, \theta_1,  v, \theta_2) \left(1+ \sum_{r = 1}^{K}N^{-r/2} b_{r}(P_0,u,v)+ N^{-(K +1)/2} R_K(P_0,u,v,N) \right), \end{multline}
where $\Pi^{T_z M}_{h_z, J_z}$ is the osculating Bargmann-Fock \szego kernel  for the tangent space $T_z M \simeq \C^m$
equipped with the complex structure $J_z$ and Hermitian metric $h_z$. Here,
\begin{itemize}

\item $ b_{r} = \sum_{\alpha=0}^{2[r/2]} \sum_{j=0}^{[3r/2]}(\psi_2)^{\alpha}Q_{r,\alpha,3r-2j}$, where $Q_{r,\alpha ,d}$ is homogeneous of
degree $d$ and
\begin{equation}
\psi_2(u,v) = u \cdot\bar{v} - \half(|u|^2 + |v|^2);
\end{equation}
(in particular, $b_r$ has only even homogeneity if $r$ is even, and only odd homogeneity if $r$ is odd);

\item $\|R_K(P_0,u,v,N)\|_{\ccal^j(\{|u|\le \rho,\ |v|\le \rho\}}\le C_{K,j,\rho}$ for $j\ge 0,\,\rho>0$ and $C_{K,j,\rho}$ is independent of the point $P_0$ and choice of coordinates.
\end{itemize}
\end{theo}

\section{Proof of logarithmic decay of variances (\autoref{theo:LOGQE})}\label{sec:LOGQE}

The variance estimate is similar to the ones given in \cite{ShZ99,Sc06,Sc08,HR16,H15}. A key ingredient is Egorov's theorem in the \kahler setting, whose proof is deferred to \autoref{sec:EGOROV}. Let $\pi \colon X \rightarrow M$ be the natural projection from the unit co-disk bundle to the base manifold. A function $f \in C^\infty(M)$ pulls back $F := \pi^* f$ to a function on $X$ that is constant along the fibers $X \rightarrow M$. Recall also that $\tilde{\chi} \colon X \rightarrow X$ is the contact lift of a symplectic diffeomorphism $\chi \colon M \rightarrow M$ for which the exponential growth estimate \eqref{eqn:EXPGROWTH} and the polynomial decay of correlations \eqref{eqn:DECAYOFCORR} apply.

\begin{prop}[Egorov's theorem with remainder]\label{prop:EGOROV}
Let $M_F$ denote multiplication by a smooth function $F := \pi^* f \in C^\infty(M)$ that is the lift of some $f \in C^\infty(M)$. Let $T \in \mathbb{Z}$ be an integer. Then
\begin{equation}
U_{\chi, N}^T (\Pi_{h^N}   M_F
 \Pi_{h^N})
 (U_{\chi, N}^*)^T= \Pi_{h^N} M_{F \circ \tilde{\chi}^T} \Pi_{h^N} + R_N^{T},
\end{equation}
where $F \circ \tilde{\chi}^T$ denotes the $T$-fold composition of $F$ with $\tilde{\chi}$, and $R_N^{T}$ is a Toeplitz operator with
\begin{equation}
\frac{1}{d_N} \Tr[(R_N^T)^*R_N^T] = \ocal\bigg(\frac{T^2}{N} \| F \|_{C^2}^2 e^{2\delta_0|T|}\bigg).
\end{equation}
In particular, at the level of matrix elements one has
\begin{equation}
\left\langle U_{\chi, N}^T \Pi_{h^N}   M_F \Pi_{h^N} (U_{\chi, N}^*)^T s_j^N , s_j^N \right\rangle = \left\langle \Pi_{h^N} M_{F \circ \tilde{\chi}^T} \Pi_{h^N} s_j^N, s_j^N \right\rangle + \ocal\bigg(\frac{T^2}{N} \| F \|_{C^2}^2 e^{2\delta_0|T|}\bigg).
\end{equation}
\end{prop}

Taking \autoref{prop:EGOROV} for granted, we proceed to prove \autoref{theo:LOGQE}. We write each integral in the Ces\`aro sum \eqref{eqn:VARIANCE} as a matrix element:
\begin{equation}\label{eqn:MATELT}
\dashint_M f(z) \|s^N_j\|_{h^N}^2\,dV = \langle \Pi_{h^N} M_F
 \Pi_{h^N} s^N_j, s^N_j \rangle.
\end{equation}
It is convenient to introduce shorthands for the time-averages:
 \begin{equation} \label{eqn:TAVE}
\begin{dcases}
[ \Pi_{h^N}  M_F
 \Pi_{h^N} ]_T : =    \frac{1}{2T+1} \sum_{n =-T}^T U_{\chi, N}^n (\Pi_{h^N}   M_F
 \Pi_{h^N})
 U_{\chi, N}^{* n},\\
[F]_{T} :=   \frac{1}{2T+1} \sum_{n =-T}^T F \circ \tilde{\chi}^n,\\
[ M_f ]_{T} := M_{[f]_{T}}.
 \end{dcases}
\end{equation}
Since $s_j^N$ are eigensections of $U_{\chi,N}$, we may replace $\Pi_{h^N} M_F \Pi_{h^N}$ in \eqref{eqn:MATELT} by its time average defined in \eqref{eqn:TAVE}:
 \begin{equation} \label{eqn:SUBSTITUTION}
 \int_M f(z) \|s^N_j\|_{h^N}^2\,dV = \left \langle [ \Pi_{h^N}   M_F
 \Pi_{h^N} ]_T s^N_j, s^N_j \right\rangle. 
\end{equation}

\autoref{prop:EGOROV}, that is Egorov's theorem, gives
\begin{equation}\label{EGTN} 
[\Pi_{h^N}   M_F
 \Pi_{h^N} ]_{T} =   \Pi_{h^N} [M_F]_{T}
 \Pi_{h^N} +  R^{(T)}_N, \end{equation}
with the remainder term satisfying the error estimate
\begin{equation} \label{RTN}
\frac{1}{d_N} \Tr[ (R^{(T)}_N)^* R^{(T)}_N] =  \ocal\bigg(\frac{T^2\|F\|^2_{C^2}e^{2\delta_0|T|}}{N}\bigg).
\end{equation}
Here the exponential growth condition \eqref{eqn:EXPGROWTH} on $\chi$ is used.

By substituting \eqref{EGTN} into \eqref{eqn:SUBSTITUTION}, the quantum variance \eqref{eqn:VARIANCE} can be rewritten as
\begin{align}
\vcal_N(f) &=   \frac{1}{d_N} \sum_{j=1}^{d_N}\left\lvert 
\left\langle [ M_F ]_{T} s^N_j, s^N_j \right\rangle + \langle R^{(T)}_N s^N_j, s^N_j \rangle -\dashint_M f\,dV
 \right \rvert^2\\
&\le \frac{2}{d_N}  \sum_{j=1}^{d_N}\left\lvert \left\langle [ M_F ]_{T} s^N_j, s^N_j \right\rangle  -\dashint_M f\,dV \right\rvert^2 + \frac{2}{d_N} \sum_{j=1}^{d_N}\left\lvert\langle R^{(T)}_N s^N_j, s^N_j \rangle \right\rvert^2.
\end{align}

Applying the Cauchy-Schwarz inequality to the first term and the error estimate \eqref{RTN} to the second term, we find 
\begin{align}\notag
\vcal_N(f) &\le \frac{2}{d_N}  \sum_{j=1}^{d_N} \dashint_M \left\lvert [f]_{T}\|s_j^N\|_{h^N}^2 - \dashint_M f\,dV \right\rvert^2 dV + \ocal\bigg(\frac{T^2\|F\|_{C^2}^2e^{2\delta_0|T|}}{N}\bigg)\\ \notag
& \le \frac{2}{d_N} \sum_{j=1}^{d_N} \dashint_M \left\lvert [f]_{T} - \dashint_M f \,dV \right\rvert^2 \|s_j^N\|_{h^N}^2\,dV+ \ocal\bigg(\frac{T^2\|F\|_{C^2}^2e^{2\delta_0|T|}}{N}\bigg)\\
& = \frac{2}{d_N} \dashint_M \left\lvert [f]_{T} - \dashint_M f \,dV \right\rvert^2 \Pi_{h^N}(z,z) \, dV + \ocal\bigg(\frac{T^2\|F\|_{C^2}^2e^{2\delta_0|T|}}{N}\bigg). \label{eqn:plugin}
\end{align}
Recall (cf.\ \cite{Z98,ShZ99}) the pointwise expansion
 for the Bergman kernel along the diagonal:
$$ \Pi_{h^N}(z,z)  =  a_0 N^m +
a_1(z) N^{m-1} + a_2(z) N^{m-2} + \cdots,$$
where the coefficients $a_j(z)$ are invariant polynomials in derivatives of the metric $h$, and where the leading order coefficient is a constant equal to $a_0 = c_1(L)^m/m!$.
Combining the Bergman kernel expansion with \eqref{eqn:plugin} yields
 \begin{equation}
\vcal_N(f)  \leq \bigg(\frac{2c_1(L)^m}{m!} + \ocal\bigg(\frac{1}{N}\bigg)\bigg)\bigg( \int_M
 \left\lvert \,  [f]_{T}  -   \dashint_M f\,dV\right\rvert^2 dV \bigg) + \ocal\bigg(\frac{T^2\|F\|^2_{C^2}e^{2\delta_0|T|}}{N}\bigg).
\end{equation}
Set
\begin{equation}
T = T(N) = \frac{1}{4\delta_0}\lvert \log N \rvert,
\end{equation}
then, thanks to the decay of correlations assumption \eqref{eqn:DECAYOFCORR}, we get (for all $0 < \beta < 1$)
 \begin{equation}
\vcal_N(f)  = \ocal\left(\frac{\|f\|_{C^{0,\beta}}^2}{\log N }\right)  +\ocal\left(\frac{\|f\|_{C^2}^2\lvert\log N\rvert^2}{N^{\frac{1}{2}}}\right) + \ocal\left(\frac{\|f\|_{C^{0,\beta}}^2}{N\log N }\right).
\end{equation}
(Note $\|F\|_{C^2} = \|f\|_{C^2}$ by definition of $F = \pi^* f$.)
This completes the proof of \autoref{theo:LOGQE}.

\section{Proofs of log-scale mass equidistribution (\autoref{cor:CESARODENSITYONE} and \autoref{theo:LOGMASS})}

\subsection{Proof of \autoref{cor:CESARODENSITYONE}}

We begin by defining constants $\kappa_1, \kappa_2$ that will appear in the proof. Let $\kappa_1$ be any constant satisfying
\begin{equation}\label{eqn:KAPPA1}
0 < \kappa_1 < 1 - 4m\gamma.
\end{equation}
It follows that 
\begin{equation}
\kappa_1 \le 1 - 4\gamma(m + \beta) \quad \text{for some $0 < \beta < 1$},
\end{equation}
whence
\begin{equation}\label{eqn:VOLUMEFACTOR}
\lvert \log N \rvert^{4\gamma\beta - 1} \le \lvert \log N \rvert^{-4m\gamma - \kappa_1}.
\end{equation} 
We also let $\kappa_2$ be any constant satisfying
\begin{equation}\label{eqn:KAPPA2}
0 < \kappa_2 < \frac{\kappa_1}{2}.
\end{equation}

Now fix $z_0 \in M$. Define symbols $\rho_N \in C^\infty_0(B(z_0,1 + 2\lvert \log N \rvert^{-\frac{\kappa_2}{\beta + 1}}, [0,1])$ by
\begin{equation}
\rho_{N}(z) := 
\begin{dcases}
1 & \text{for $z \in B(z_0,1 + \lvert \log N \rvert^{-\frac{\kappa_2}{\beta + 1}})$,}\\
0 & \text{for $z \notin B(z_0,1 + 2\lvert \log N \rvert^{-\frac{\kappa_2}{\beta + 1}})$.}
\end{dcases}
\end{equation}
Note that the support of $\rho_N$ depends on $N$. We perform a further rescaling
\begin{equation}\label{eqn:RHON}
(D_{\epsilon_N}^{-1})^*\rho_N(z) = \rho_N(\epsilon_N^{-1} z).
\end{equation}
The statement of \eqref{cor:LOGQE} (which follows easily from Theorem~\autoref{theo:LOGQE} as discussed in \autoref{sec:QEintro})  with $f_{z_0,\epsilon_N}$ replaced by $\rho_N(\epsilon_N^{-1}z)$ becomes
\begin{align}
\frac{1}{d_N} \sum_{j=1}^{d_N} \left\lvert \big\langle M_{(D_{\epsilon_N}^{-1})^* \rho_N} s_j^N, s_j^N \big\rangle  - \dashint_M \rho_N(\epsilon_N^{-1}z)\,dV \right\rvert^2 &= \ocal(\|\rho_N\|_{C^{0,\beta}}^2 \lvert \log N \rvert^{4\gamma\beta -1 })\\
&\le \ocal(\|\rho_N\|_{C^{0,\beta}}^2\lvert \log N \rvert^{-4m\gamma}\lvert \log N \rvert^{-\kappa_1}) \label{eqn:SMALLSCALEQE2}
\end{align}
for any $\kappa_1$ satisfying \eqref{eqn:KAPPA1}. In the last line we used \eqref{eqn:VOLUMEFACTOR}.

Now apply Markov's inequality
$\mathbb{P}(X \geq a) \leq a^{-1} {\mathbb{E}} X$. We view each term of the sum on the left-hand side of \eqref{eqn:SMALLSCALEQE2} as a random variable indexed by $(N,j)$. The probability measure is the normalized counting measure on the indices $\{0 \le j \le d_N\}$. Finally take $a$ to equal $\lvert \log N \rvert^{\epsilon}$ (for some small $\epsilon > 0$) times the right side of \eqref{eqn:SMALLSCALEQE2}. It follows that  for any constant $\kappa_2$ satisfying \eqref{eqn:KAPPA2}
there exists a full density subsequence $\Gamma_{z_0}' \subset \{(N,j)\}$ such that the corresponding eigensections satisfy
\begin{equation}\label{eqn:DENSITYONE}
\left\lvert \int_{B(z_0,2)} \rho_N(\epsilon_N^{-1}z) \| s_j^N\|_{h^N}^2 - \frac{1}{\Vol(M)} \int_{B(z_0,2)} \rho_N(\epsilon_N^{-1}z)\, \right\rvert
\le C \|\rho_N\|_{C^{0,\beta}}\lvert \log N \rvert^{-2m\gamma}\lvert \log N \rvert^{-\kappa_2}
\end{equation}
for $(N,j) \in \Gamma'_{z_0}$.
In other words, almost all the terms in the averaged sum \eqref{eqn:DENSITYONE} each satisfies the slightly worse than  the average  upper bound $C \|\rho_N\|_{C^{0,\beta}}\lvert \log N \rvert^{-2m\gamma}\lvert \log N \rvert^{-\kappa_2}$.

We then have
\begin{align}
\int_{B(z_0,\epsilon_N)} \| s_j^N \|_{h^N}^2\,dV &\le \int_{B(z_0,2)} \rho_N(\epsilon_N^{-1}z) \| s_j^N \|_{h^N}^2\,dV\\
&\le \frac{1}{\Vol(M)} \int_{B(z_0,2)} \rho_N(\epsilon_N^{-1}z)\,dV + C\|\rho_N\|_{C^{0,\beta}}\lvert \log N \rvert^{-2m\gamma}\lvert \log N \rvert^{-\kappa_2}\\ \label{eqn:ESTIMATE}
&\le \frac{\Vol(B(z_0, \epsilon_N))}{\Vol(M)} + C\left( \lvert \log N \rvert^{-2m\gamma-\frac{\kappa_2}{\beta +1}}+ \|\rho_N\|_{C^{0,\beta}}\lvert \log N \rvert^{-2m\gamma-\kappa_2}\right).
\end{align}
The first inequality follows from the definition \eqref{eqn:RHON} of $\rho_N$. The second inequality follows from the estimate \eqref{eqn:DENSITYONE}. 
The third inequality follows from the support condition of \eqref{eqn:RHON} and from the  volume of spherical shells (the ``thickness'' of the shell being $2\lvert \log N \rvert^{-\frac{\kappa_2}{\beta +1}} $):
\begin{align}
\int_{B(z_0,2)} \rho_N(\epsilon_N^{-1}z)\,dV &= \int_{B(z_0,1+2\lvert \log N \rvert^{-\frac{\kappa_2}{\beta +1}})\setminus B(z_0,1)}\rho_N(\epsilon_N^{-1}z) \,dV + \int_{B(z_0,1)}\rho_N(\epsilon_N^{-1}z)\,dV\\
&\le \epsilon_N^{2m}\int_{B(z_0,1+2\lvert \log N \rvert^{-\frac{\kappa_2}{\beta +1}})\setminus B(z_0,1)} dV + \int_{B(z_0,\epsilon_N)} dV\\
&\le C \epsilon_N^{2m}\lvert \log N \rvert^{-\frac{\kappa_2}{\beta +1}} + \Vol(B(z_0,\epsilon_N)),
\end{align}
where $C$ depends only on $(M,\omega)$ and the choice of $\rho$. 

Note that $\|\rho_N\|_{C^{0,\beta}} \le C (\lvert \log N \rvert^{\frac{\kappa_2}{\beta+1}})^{-\beta}$, 
which gives
\begin{align}
\int_{B(z_0,\epsilon_N)} \| s_j^N \|_{h^N}^2\,dV &\le \frac{\Vol(B(z_0, \epsilon_N))}{\Vol(M)} + C \lvert \log N \rvert^{-2m\gamma}\bigg(\lvert \log N \rvert^{-\frac{\beta\kappa_2}{\beta + 1}} + \lvert \log N \rvert^{-\frac{\beta\kappa_2}{\beta+1}}\bigg)\\
&= \frac{\Vol(B(z_0, \epsilon_N))}{\Vol(M)} + \smallO(\lvert \log N \rvert^{-2m\gamma}). \label{eqn:LE}
\end{align}
(From \eqref{eqn:KAPPA1} and \eqref{eqn:KAPPA2} of how $\kappa_1,\kappa_2$ are defined, we have $0 < \beta\kappa_2/(\beta+1) < 1$.)

A similar argument using appropriately chosen $\tilde{\rho}_N$ of the form $\tilde{\rho}_N(z) = \rho_N(3z)$ gives the opposite inequality
\begin{equation}\label{eqn:GE}
\int_{B(z_0,\epsilon_N)} \| s_j^N \|_{h^N}^2\,dV \ge \frac{\Vol(B(z_0, \epsilon_N))}{\Vol(M)} + \smallO(\lvert \log N \rvert^{-2m\gamma})
\end{equation}
for a full density subsequence $\Gamma_{z_0}''$ of eigensections. The intersection $\Gamma_{z_0}' \cap \Gamma_{z_0}'' =: \Gamma_{z_0}$ indexes a full density subsequence of eigensections for which \eqref{eqn:LE} and \eqref{eqn:GE} hold simultaneously. This completes the proof of \autoref{cor:CESARODENSITYONE}.

\subsection{Proof of \autoref{theo:LOGMASS}}

Note that one must first fix a single base point $z_0 \in M$ for the asymptotic statement of \autoref{cor:CESARODENSITYONE} to hold. To move towards global statements that hold for all $z \in M$ simultaneously, we introduce the concept of a \emph{log-good cover}, for which we have uniform estimates on each element (i.e., a \kahler ball) of the cover. The existence of a cover satisfying the following conditions is proved in \cite{H15}.

\begin{defin}
Let $\epsilon_N = \lvert \log N \rvert^{-\gamma}$ for any fixed $0 < \gamma < (6m)^{-1}$ as before. A log-good cover $\ucal_N$ is a cover of $M$ by geodesic balls $\{B(z_{N,\alpha},\epsilon_N)\}_{\alpha = 1}^{R(\epsilon_N)}$ with the following properties:
\begin{itemize}
\item The number $R(\epsilon_N)$ of balls in the cover is bounded above
\begin{equation}
R(\epsilon_N) \le c_1 \epsilon_N^{-2m} \qquad (\dim_\R M = 2m)
\end{equation}
by some constant \textup{(}independent of $N$\textup{)} multiple of $\epsilon_N^{-2m}$.

\item An arbitrary ball $B(p,\epsilon_N) \subset M$ is covered by at most $c_2$ \textup{(}independent of $N$\textup{)} number of balls from the cover. 

\item An arbitrary ball $B(p,\epsilon_N) \subset M$ contains at least one of the shrunken balls $B(z_{N,\alpha},\frac{\epsilon_N}{3})$.
\end{itemize}
\end{defin}

We now proceed with the proof of \autoref{theo:LOGMASS}, suppressing the prime notation on $\gamma$ and $\epsilon_N$.
Let $0 < \gamma < (6m)^{-1}$ be given and set $\epsilon_N = \lvert \log N \rvert^{-\gamma}$. For each $N$, fix a log-good cover $\ucal_N$ as defined above. As before, let $0 \le f_{z_{N,\alpha}} \le 1$ be a smooth cut-off function that is equal to 1 on $B(z_{N,\alpha},1)$, and vanishes outside $B(z_{N,\alpha},2)$. Let $f_{z_{N,\alpha},\epsilon_N} = f_{z_{N,\alpha}}(\epsilon_Nz)$. (This is a slight abuse of notation, where we mean balls in \kahler normal coordinate charts centered at $z_{N,\alpha}$.) In what follows, $\kappa_3 > 0$ is a parameter independent of $N,j$, to be chosen later.

The extraction argument uses Markov's inequality $\mathbb{P}(X \geq a) \leq a^{-1}{\mathbb{E}} X$. To this end, for each $1 \le j \le d_N$ and $1 \le \alpha \le R(\epsilon_N)$ set
\begin{equation}
X_{N,j,\alpha} := \left\lvert\, \int_M f_{z_{N,\alpha},\epsilon_N}\|s^N_j\|_{h^N}^2\,dV  - \dashint_M f_{z_{N,\alpha},\epsilon_N}\,dV \,\right\rvert^2.
\end{equation}
We view $X_{N,j,\alpha}$ as a random variable with respect to the normalized counting measure on the set of indices $1 \le j \le d_N$. Thanks to \autoref{cor:LOGQE} and \eqref{eqn:VOLUMEFACTOR}, its expected value is
\begin{equation}
\mathbb{E}X_{N,j,\alpha} = \ocal(\lvert \log N \rvert^{-(1-2\gamma\beta)}) = \ocal(\lvert \log N \rvert^{-(4m\gamma + \kappa_1)}) \quad \text{for any $\kappa_1$ satisfying \eqref{eqn:KAPPA1}}.
\end{equation}
(The error is uniform in $z_{N,\alpha}$.) In particular, we may choose $\kappa_1$ to equal
\begin{equation}\label{eqn:KAPPA1DEFINED}
0 < \kappa_1 := 1- 4m(\gamma + \beta)  < 1 \quad \text{for some $0 < \beta < \frac{1-6m\gamma}{4m} < 1$}.
\end{equation}
It follows from an application of Markov's inequality with $X = X_{N,j,\alpha}$; with the normalized counting measure on $\{1, \dots, d_N\}$; and with $a = \lvert \log N \rvert^{-(4m\gamma-\kappa_3)}$, that the `exceptional sets' 
\begin{equation}
\Lambda_{\alpha}(N) := \bigg\{ j = 1, \dotsc, d_N \colon 
\left\lvert\, \int_M f_{z_{N,\alpha},\epsilon_N}\|s^N_j\|_{h^N}^2\,dV  - \dashint_M f_{z_{N,\alpha},\epsilon_N}\,dV\right\rvert^2 \ge \lvert \log N \rvert^{-4m\gamma-\kappa_3}\bigg\}
\end{equation}
satisfy
\begin{equation}
\frac{\# \Lambda_{\alpha}(N)}{d_N} \leq C \lvert \log N \rvert^{4m\gamma-\kappa_3}\lvert \log N \rvert^{-(4m\gamma+\kappa_1)} = C \lvert \log N \rvert^{-(1-4m(\gamma+\beta) - \kappa_3)}.
\end{equation}

Now define `generic sets'  
\begin{equation} \label{eqn:GENERICSET}
\Sigma_\alpha(N) : = \{j \colon 1 \leq j \leq d_N\} \setminus \Lambda_{\alpha}(N) \quad \text{and} \quad \Sigma(N) := \bigcap_{\alpha \colon B(z_{N,\alpha},\epsilon_N) \in \ucal_N} \Sigma_{\alpha}(N).
\end{equation} 
The number of elements in the cover $\ucal_N$ is of order $\epsilon_N^{-2m} = \lvert \log N \rvert^{2m\gamma}$, whence
\begin{align}
\frac{\# \Sigma(N)}{d_N} &\geq 1 - \sum_{\alpha} \frac{\#\Lambda_{\alpha}(N)}{d_N} \\
&\ge 1 - C|\log N|^{2m \gamma} \lvert\log N \rvert^{-(1-4m(\gamma+\beta) - \kappa_3)}\\
& = 1 - C \lvert \log N \rvert^{-(1-6m\gamma-4m\beta-\kappa_3)}\\
& \rightarrow 1 \quad \text{by choosing $\beta,\kappa_3 > 0$ sufficiently small.}\label{eqn:GOTOZERO}
\end{align}
Indeed, by choice \eqref{eqn:KAPPA1DEFINED} of $\beta$, we have $1-6m\gamma-4m\beta > 0$, so
$\kappa_3$ can always be chosen to ensure \eqref{eqn:GOTOZERO} holds. This is
analogous to the estimate in \cite{HR16} preceding Lemma~3.1 or in \cite[p.3263]{H15}.

The construction of indexing sets $\Sigma(N)$ yields a full density subsequence
\begin{equation}
\Sigma := \bigcup_{N \ge 1} \Sigma(N)
\end{equation}
such that, for every $B(z_\alpha,\epsilon_N) \in \ucal_N$, we have 
\begin{align}
\int_{B(z_{N,\alpha},\epsilon_N)} \|s_j^N\|_{h^N}^2\,dV &\le \int_{B(0,2)} f_{z_{N,\alpha},\epsilon_N} \|s_j^N\|_{h^N}^2\,dV\\
& \le \frac{1}{\Vol(M)} \int_{B(0,2)} f_{z_{N,\alpha},\epsilon_N}\,dV + C\lvert \log N \rvert^{-(2m\gamma+\kappa_3/2)}\\
& \le \frac{\Vol(B(z_{N,\alpha},2\epsilon_N))}{\Vol(M)} + \smallO(\lvert \log N \rvert^{-2m\gamma})\\
&\le C \Vol(B(z_{N,\alpha},\epsilon_N))
\end{align}
simultaneously for all $\alpha = 1, \dotsc, R(\epsilon_N)$ as $\Sigma \ni (N,j) \rightarrow \infty$. The constant $C$ is independent of $\alpha$.

Now let $p \in M$ be arbitrary. By construction, the ball $B(p,\epsilon_N)$ is contained in at most $c_2$ number (independent of $N$)  of elements of the log-good cover $\ucal_N$. Thus, 
\begin{equation}
\int_{B(p,\epsilon_N)} \|s_j^N\|_{h^N}^2\,dV  \le \sum_{i=1}^{c_2}\frac{1}{\Vol(M)} \int_{B(0,2)} f_{z_{N,\alpha_i},\epsilon_N}\,dV + \smallO(\lvert \log N \rvert^{-2m\gamma})\le C\Vol(B(p,\epsilon_N))
\end{equation}
for every $p \in M$ as $\Sigma \ni (N,j) \rightarrow \infty$. The constant $C$ is independent of $p$. This is the statement of the volume upper bound.

It remains to repeat the same construction by dilating the symbol $0 \le g_{z_\alpha} \le 1$ that is a smooth cut-off function supported in $B(z_\alpha,1/3)$ and equals to 1 in $B(0, 1/6)$. There exists a full density subsequence $\Sigma'$ such that
\begin{align}
\int_{B(z_{N,\alpha}, \epsilon_N/3)} \|s_j^N\|_{h^N}^2\,dV & \ge \int_{B(z_\alpha,1/3)} g_{z_\alpha,\epsilon_N} \|s_j^N\|_{h^N}^2\,dV\\
&\ge \frac{1}{\Vol(M)}\int_{B(z_{N,\alpha},1/3)} g_{z_\alpha,\epsilon_N/3} \,dV - C\lvert \log N \rvert^{-(2m\gamma+\kappa_3/2)}\\
& \ge \frac{\Vol(B(z_{N,\alpha},\epsilon_N/6))}{\Vol(M)} - \smallO(\lvert \log N \rvert^{-2m\gamma})\\
&\ge c \Vol(B(z_{N,\alpha},\epsilon_N))
\end{align}
simultaneously for all $\alpha = 1, \dotsc, R(\epsilon_N)$ as $\Sigma \ni (N,j) \rightarrow \infty$. Now let $p \in M$ be arbitrary. Every ball $B(p,\epsilon_N)$ contains at least one element $B(z_{N,\alpha},\epsilon_N/3) \in \ucal_N$ of the log-good cover, whence
\begin{equation}
\int_{B(p,\epsilon_N)} \|s_j^N\|_{h^N}^2\,dV \ge c \Vol(B(p,\epsilon_N))
\end{equation}
for every $p \in M$ as $\Sigma \ni (N,j) \rightarrow \infty$. This is the statement of the volume lower bound.

The intersection $\Gamma = \Sigma \cap \Sigma'$ is again a full density subsequence. By construction, the eigensections indexed by $\Gamma$ satisfy the two-sided bound: for all $p \in M$,
\begin{equation}
c \Vol(B(p,\epsilon_N)) \le \int_{B(p,\epsilon_N)} \|s_j^N\|_{h^N}^2\,dV \le C \Vol(B(p,\epsilon_N)) \quad \text{as $\Gamma \ni (N,j) \rightarrow \infty$}. 
\end{equation}
This completes the proof of \autoref{theo:LOGMASS}.

\subsection{Proof of log-scale equidistribution of zeros (\autoref{theo:LOGZERO})}
Let $0 < \gamma < (6m)^{-1}$ from the statement of \autoref{theo:LOGZERO} be given. 
We distinguish two logarithmic scales by fixing another parameter $\gamma'$:
\begin{equation}
0 < \gamma < \gamma' < \frac{1}{6m} \quad \text{so that} \quad \lvert \log N \rvert^{-\gamma'} = \epsilon_N' < \epsilon_N = \lvert \log N \rvert^{-\gamma}.
\end{equation}
Let $\Gamma$ be the full density subsequence corresponding to scale $\epsilon'$ as guaranteed by \autoref{theo:LOGMASS}. We show that the same $\Gamma$ satisfies the statement of \autoref{theo:LOGZERO} at the scale $\epsilon_N > \epsilon_N'$.

In the notation of \autoref{CG}, relative to a local frame we write the eigensections locally as
\begin{equation}
s_j^{N} = f_j^{(N)}e_L^N, \quad \text{$f_j^{(N)}$ a local holomorphic function.}
\end{equation}
The Poincar\'e-Lelong formula \eqref{PL} reduces the growth rate of zeros to the growth rate of the local plurisubharmonic function
$N^{-1} \log |f_j^{(N)}|^2$ or  to the global quasi-plurisubharmonic function\footnote{`quasi' means p.s.h. up to a fixed continuous term, here the potential $\log g$ where $g(z) := \|e_L(z)\|_h^2$.} $u^{(N)}_j(z) = N^{-1} \log \|s_j^N (z)\|_{h^N}^2$. Fix $p \in M$ and  consider the dilated
function
\begin{equation}\label{DILATEDu} 
u_j^{(N)}(z) := \frac{1}{N} \log \|s_j^{N}(\epsilon_{N} z)\|_{h^N}^2 = D_{\epsilon_{N}}^{p*}\bigg[ \frac{1}{N} \log \|s_j^{N}( z)\|_{h^N}^2\bigg] \quad\text{on $B(p, 1)$},
\end{equation}
where  $D^p_{\epsilon_{N}}$ is the local dilation defined by \eqref{Dep}
 in \kahler normal coordinates centered at $p = 0$. Since $D_{\epsilon_{N}}^p$ is a local holomorphic map, \eqref{DILATEDu} remains
quasi-plurisubharmonic. We state a key lemma:

\begin{lem}\label{W*}
Let   $\Gamma$ be the subsequence of density one  for the finer scale $\epsilon_N'$ of  \autoref{theo:LOGMASS}. For $(N,j) \in \Gamma $, the logarithmically dilated potential
\eqref{DILATEDu} satisfies
\begin{equation}
\| u_j^{(N)} \|_{L^1(B(p,1))} = \smallO(\epsilon_N^2),
\end{equation}
where the remainder is at a coarser scale $\epsilon_N$.
\end{lem}

\begin{rem}
We emphasize that we are assuming the eigensections indexed by $\Gamma$ satisfy
\begin{equation}\label{CONTRA}
 C_1 \frac{\Vol(B(p, \epsilon_N'))}{\Vol(M)} \le \int_{B(p, \epsilon_N')} \| s_j^N\|_{h^N}^2 \, dV \le C_2 \frac{\Vol(B(p, \epsilon_N'))}{\Vol(M)}
\end{equation}
and then inverse dilating $B(p,\epsilon_N)$ to $B(p,1)$, so that  any ball $B(q, \epsilon_N')
\subset B(p, \epsilon)$ gets inverse dilated to (slightly deformed) by $(D^{p}_{\epsilon_N})^{ -1}$  to (slightly deformed) balls of radius $\epsilon_N^{-1} \epsilon_N' \simeq |\log N|^{-\gamma' + \gamma}$ in $B(p, 1)$. 
\end{rem}

Let's assume \autoref{W*} for now and proceed to finish the proof of \autoref{theo:LOGZERO}. Using the Poincar\'e-Lelong formula and the fact that the holomorphic rescaling $D_\epsilon^p$ commutes with $\ddbar$, we obtain
\begin{equation}\label{LAST}
\frac{1}{N} D_{\epsilon_{N}}^{p*} \big[Z_{s_j^N}\big] = \frac{\sqrt{-1}}{2\pi N} \ddbar \log |f^{(N)}_j(\epsilon_N z)|^2 = \frac{\sqrt{-1}}{2\pi N}\ddbar \log \| s^{N}_j(\epsilon_N z) \|^2_{h^N}  + D_{\epsilon_N}^{p*}\omega.
\end{equation} 
For every test form $\eta \in \dcal^{m-1,m-1}(B(p,1))$ and $\Gamma \ni (N,j) \rightarrow \infty$, integration by parts
and \autoref{W*} give
\begin{align}
\int_{B(p,1)}  \left(\eta \wedge \frac{1}{N} D_{\epsilon_{N}}^{p*} \big[Z_{s_j^N}\big] \right) &= \int_{B(p,1)}\eta  \wedge D_{\epsilon_{N}}^{p*}\omega + \int_{B(p,1)} \frac{\sqrt{-1}}{2\pi N}\log \| s^{N}_j(\epsilon_N z) \|^2_{h^N} \ddbar\eta(z) \\
& = \int_{B(p,1)}\eta  \wedge D_{\epsilon_{N}}^{p*}\omega + \smallO(\epsilon_N^2).\label{eqn:BYPARTS}
\end{align}
Locally at $p = 0$, the \kahler potential can be written as $\phi(z) = \lvert z \rvert^2 + \ocal(\lvert z \rvert^4)$, so
\begin{equation}\label{eqn:LOCALPHI}
D_{\epsilon_N}^{p*}\omega = \frac{\sqrt{-1}}{2\pi} D_{\epsilon_N}^{p*} \ddbar \phi = \epsilon_N^2 \frac{\sqrt{-1}}{2\pi}\ddbar \lvert z \rvert^2 + \ocal(\epsilon_N^4) = \epsilon_N^2 \omega_0^p + \ocal(\epsilon_N^4),
\end{equation}
with $\omega_0^p$ the flat \kahler form. Combining \eqref{eqn:BYPARTS} and \eqref{eqn:LOCALPHI} (and dividing by $\epsilon_N^2$) yields
\begin{equation}
\int_{B(p,1)}  \left(\eta \wedge \frac{1}{N\epsilon_N^2} D_{\epsilon_{N}}^{p*} \big[Z_{s_j^N}\big] \right) = \int_{B(p,1)}\eta  \wedge \omega_0^p + \smallO(1) \quad \text{as $\Gamma \ni (N,j) \rightarrow \infty$},
\end{equation}
which is equivalent to the statement of \autoref{theo:LOGZERO}.

\begin{proof}[Proof of \autoref{W*}]
The 
argument is similar to the one in \cite{ShZ99} except for the dilation
of the plurisubharmonic functions. 
The  log-scale quantum ergodicity successfully replaces
 unscaled quantum ergodicity in the key step of the argument due to the
 fact that the local dilation is holomorphic. But we need to use two logarithmic scales and for later applications we need the remainder estimate. 
 
Let $N_0$ be sufficiently large so that for all $N \ge N_0$, $e_L$ is a local frame for $L$ over an open subset $U$ containing $\overline{B(p,1)} $ and 
$e_L^{N}$ is the corresponding frame for $L^{N}$. Since  $g(z) = \|e_L(z)\|_h^2$, we have
\begin{equation}
\|e_L^{N}(z)\|_{h^N}^2 =  g^N  \quad \text{and} \quad \|s_j^{N}(\epsilon_{N} z)\|_{h^N}^2 = |f_j^{N}(\epsilon_{N} z)|^2 g^{N}(\epsilon_N z).
\end{equation} 
We first show that $\|u_j^{(N)}\|_{L^1} \rightarrow 0$, and then indicate how the argument can be adapted to yield the $\smallO(\epsilon_N^{2})$ improvement.

Observe that any $L^2$-normalized section
satisfies
 \begin{equation} \label{SNUB} \| s^N
(z)\|_{h^N}^2 \leq \Pi_{h^N}(z,z)= \left(\frac{c_1(L)^m}{m!} + O\bigg(\frac{1}{N}\bigg)\right)N^m.\end{equation} Hence
$\| s^N (z)\|_{h^N} \leq C N^{m/2}$ for some $C <\infty$ and taking the
logarithm gives
\begin{itemize}
\item[(i)] The functions $u^{(N)}$ are
uniformly bounded above on $M$;
\item[(ii)] $\limsup_{N \rightarrow \infty} u_N \leq 0$.
\end{itemize} 
Now consider the plurisubharmonic function 
\begin{equation}
v_j^{(N)}(z) := \frac{1}{N}\log \lvert f_j^{(N)}(\epsilon_{N} z) \rvert^2 = u_j^{(N)}(z) - \log g (\epsilon_{N} z) \in \operatorname{PSH}(B(p,1)).
\end{equation}
It is clear that $v_j^{(N)}$ are uniformly upper bounded. A standard result on plurisubharmonic functions (see
\cite[Theorem~4.1.9]{H90}) then implies a subsequence $v_j^{(N_k)}$ either
converges uniformly to $-\infty$ on $B(p,1)$ or else has a subsequence that is
convergent in $L^1_\mathrm{loc}(B(p,1))$.

Let us  rule out the first possibility. If it
occurred,  there would exist $K > 0$ such that
\begin{equation}\label{firstposs}
\frac{1}{N_k} \log \| s_j^{N_k}(\epsilon_{N_k} z)\|_{h^{N_k}}^2 \leq -
1 \iff \| s_j^{N_k}(\epsilon_{N_k} z)\|_{h^{N_k}}^2 \leq e^{- N_k} \quad \text{on $B(p,1)$ for all $k \ge K$}.
\end{equation}
Equivalently, the same exponential decay estimate holds on $B(p, \epsilon_{N_k})$ for the undilated sections. But this contradicts the lower bound 
of \eqref{CONTRA}.

Therefore the sequence  $v_j^{(N)}$ is pre-compact in $L^1(B(p,1))$, and  every sequence contains  a subsequence, which we continue to denote  by
$\{v_j^{(N_k)}\}$,
that converges in $L^1(B(p, 1))$ to some $v
\in L^1(B(p,1))$.  By
passing if necessary to a further subsequence, we may assume that
$\{v^{(N_k)}_j \}$ converges pointwise
almost everywhere  in $B(p,1)$ to $v$, and hence by observation (ii),
\begin{equation}
v(z) = \limsup_{(N_k,j) \rightarrow \infty} \left(u_j^{(N_k)}(z)- \log g(\epsilon_{N_k}z) \right) \leq 0 \quad \text{a.e.\ on $B(p,1)$.}
\end{equation}

Let
\begin{equation}
v^*(z):= \limsup_{w \rightarrow z} v(w) \leq  0
\end{equation} be the
upper-semicontinuous regularization of $v$. Then $v^*$ is plurisubharmonic on
$B(p,1)$ and $v^* = v$ almost everywhere.
We claim that $v^* =0$.
To this end, we use the second scale $\epsilon_N'$.
If $v^* \not= 0$, then
\begin{equation}
\|v_j^{(N_k)} + D_{\epsilon_{N_k}}^{p*}\log g\|_{L^1(B(p,1))}=\|u_j^{(N_k)}\|_{L^1(B(p,1))} \geq
\delta>0.
\end{equation}
Hence, for some $c
> 0$, the open set $U_c = \{z\in B(p,1) \colon v^*(z) <  - c\} $ is
nonempty. For sufficiently large $k$, this set contains a ball $B(q, \epsilon_{N_k}' \epsilon_{N_k}^{-1})$. By Hartogs' Lemma, there exists
a positive integer $K$ such that $v_j^{(N_k)}(z) \leq -c/2$ for $z \in
B(q, \epsilon_{N_k}' \epsilon_{N_k}^{-1})$ and $k\geq K$, that is
\begin{equation}
\|s_j^{N_k}(\epsilon_{N} z)\|_{h^{N_k}}^2 \leq e^{-c N_k/2 } \quad \text{on $B(q, \epsilon_{N_k}' \epsilon_{N_k}^{-1})$ for all $k \geq K$.}
\end{equation}
But this again  contradicts the lower bound in \autoref{theo:LOGMASS} on $B(q, \epsilon_{N_k}')$. We have therefore proved $\|u_j^{(N)}\|_{L^1(B(p,1))} = \smallO(1)$.

We now exploit the exponential decay to 
prove the sharper result $\|u_{j}^{(N)}\|_{L^1(B(p,1))} = \smallO(
\epsilon_N^2)$.  Consider the renormalized sequence
\begin{equation}
\epsilon_N^{-2} u_j^{(N)}
= \frac{1}{N \epsilon_N^2} D_{\epsilon_N}^* \log \|s^N_j(z)\|^2_{h^N}.
\end{equation}
Note that  this
 is still an upper-bounded sequence of plurisubharmonic functions because of the exact cancellation between dilating by $D_{\epsilon_N}^{p*}$ and dividing by $\epsilon_N^{2}$. 
Indeed, $\log g = |z|^2 + \ocal(|z|^4)$ as $|z| \to p = 0$ in local coordinates, so $\epsilon_N^{-2} D_{\epsilon_N}^{p*} \log g$ remains bounded.

We now run
through the previous argument again with this re-normalized sequence. 
If $\epsilon_{N_k}^{-2} v_j^{N_k} \rightarrow -\infty$ uniformly on compact subsets of $B(p,1)$, then
\begin{equation}
\frac{1}{N_k\epsilon_{N_k}^2}\|s^{N_k}_j(\epsilon_{N_k}z)|_{h^{N_k}}^2
\leq  -1 \iff \|s^{N_k}_j(\epsilon_{N_k}z)\|_{h^{N_k}}^2 \le e^{-\epsilon_{N_k}^2N_k} \quad \text{on $B(p,1)$},
\end{equation}
a contradiction to \eqref{CONTRA} as before.
The alternative (namely $\epsilon_{N_k}^{-2}v^{N_k}_j$ being pre-compact) leads to the estimate
\begin{equation}
\|s_j^{N_k}(\epsilon_{N} z)\|_{h^{N_k}}^2 \leq e^{-c\epsilon_{N_k}^2 N_k/2 } \quad \text{on $B(q, \epsilon_{N_k}' \epsilon_{N_k}^{-1})$ for all $k \geq K$,}
\end{equation}
again a contradiction. This completes the proof of \autoref{W*}.
\end{proof}

\appendix

\section{Egorov's theorem} \label{sec:EGOROV}

The purpose of this section is to prove a long time  Egorov's theorem with 
remainder as stated in \autoref{prop:EGOROV}. It is convenient to work on the contact manifold $(X,\alpha)$ by lifting $\chi$ on $M$ to the contact transformation $\tilde{\chi}$ on $X$ and viewing sections $s_j^N \in H^0(M,L^N)$ as equivariant functions $\hat{s}_j^N \in L^2(X)$ as discussed in \autoref{sec:BACKGROUND}.

We recall the setting. Let $\chi$ be a quantizable symplectic map (whose quantization $U_{\chi,N}$ is defined in \eqref{eqn:QUANTIZE1}) satisfying the exponential growth condition \eqref{eqn:EXPGROWTH} and decay of correlations condition \eqref{eqn:DECAYOFCORR}. Let
 $M_F$ denote multiplication by $F \in C^{\infty}(X)$ and $F \circ \tilde{\chi}^T$ the composition of $F$ with the $T$-fold iterate of $\tilde{\chi}$ (or $\tilde{\chi}^{-1}$, depending on the sign of $T$). \autoref{prop:EGOROV}, which is a statement on the base manifold $M$, is equivalent to the following statement on the co-circle bundle $X$.

\begin{prop}\label{prop:EGOROV2}
Let $\chi$ be a quantizable symplectic map on $M$ satisfying conditions \eqref{eqn:EXPGROWTH} and \eqref{eqn:DECAYOFCORR}. Let $\tilde{\chi}$ denote its lift to $(X,\alpha)$ as a contact transformation. Let $F \in C^\infty(X)$ and $T \in \mathbb{N}$. Then
\begin{equation}
U_{\chi, N}^T (\Pi_{h^N}   M_F \Pi_{h^N}) (U_{\chi, N}^*)^T= \Pi_{h^N} M_{F \circ \tilde{\chi}^T} \Pi_{h^N} + R_N^{(T)},
\end{equation}
where $R_N^{(T)}$ is a Toeplitz operator with
\begin{equation}
\frac{1}{d_N} \|R_N^{(T)} \|_{\mathrm{HS}}^2 = \frac{1}{d_N} \Tr[(R_N^{(T)})^*R_N^{(T)}] = \ocal_{\tilde{\chi}, F, h}\bigg(\frac{T^2}{N}\|F\|_{C^2}^2 e^{2\delta_0|T|}\bigg),
\end{equation}
where the $\ocal$ symbol depends on the metric $h$ and a fixed number of derivatives
of $\tilde{\chi}, F$ depending on the dimension.
\end{prop}

The proposition is the analogue for Toeplitz operators of the well-known estimate of the Egorov remainder, except that the remainder is stated in terms
of the normalized Hilbert-Schmidt norm rather than the operator norm.\footnote{The more difficult norm estimate of the remainder will be presented elsewhere.} The Hilbert-Schmidt norm is simpler to estimate since it is defined by a trace, and the remainder estimate is simply the standard one in the stationary phase expansion \cite{H90}.   Sharper  remainder estimates have been proved for quantizations of Hamiltonian flows on $T^*\R^n$ in Theorems~1.4
and 1.8 of  \cite{BouR02}. Subsequently, there are many articles proving related results for $T^*M$. But there do not seem to exist parallel results for Toeplitz operators in the \kahler setting, in particular for powers of a map rather than for Hamiltonian flows. In special cases such as symplectic toral automorphisms and their perturbations, Egorov's theorem with remainder have been proved (see \cite{Sc06,Sc08}) but the proofs use special properties of the metaplectic representation and do not generalize to our setting.  Egorov's theorem without estimate of the time-dependence of the remainder  may be obtained from the composition theorem for Toeplitz operators in \cite{BG81}.

\begin{rem} The strategy of the proof is to use induction on $T$. At each
stage, the remainder terms from the previous stage are left `untouched', and are estimated using that  unitary conjugations do not change Hilbert-Schmidt norms. Unlike most statements of the Egorov theorem, we only need the principal term and
a remainder of order $N^{-1}$, and we do not try to give a formula for
the lower order terms in the symbol. Thus, at the $T$th stage we only
conjugate by one power of $U_{\chi, N}$ a Toeplitz operator whose symbol 
is of the form $F \circ \tilde{\chi}^{T-1}$. This is why the resulting remainder after $T$ steps involves the $C^2$ norm of $F \circ \tilde{\chi}^{T}$ and otherwise only involves
a fixed number of derivatives of the data $\tilde{\chi}, h, F$. \end{rem}

\subsection{Reduction to \texorpdfstring{$T = 1$}{T = 1} case}

In this section we reduce the proof of \autoref{prop:EGOROV2} to the proof of the following lemma.

\begin{lem} \label{lem:SUFFICE}
Under the same assumption as \autoref{prop:EGOROV2}, we have
\begin{equation}\label{eqn:SUFFICE}
U_{\chi, N} \Pi_{h^N}   M_F \Pi_{h^N} U_{\chi, N}^* = \Pi_{h^N} M_{F \circ \tilde{\chi}} \Pi_{h^N} + R_N,
\end{equation}
where $R_N$ is a Toeplitz operator with
\begin{equation}
\frac{1}{d_N}\left\|
R_N\right\|_{\mathrm{HS}}^2 = \frac{1}{d_N} \Tr[R_N^*R_N] = \ocal_{\tilde{\chi}, F, h}\bigg(\frac{1}{N}\| F\|_{C^2}^2 e^{2\delta_0}\bigg).
\end{equation}
\end{lem}

We now indicate how \autoref{lem:SUFFICE} implies the statement of Egorov's theorem. The rest of the appendix is then devoted to proving \autoref{lem:SUFFICE}.

\begin{proof}[Proof of \autoref{prop:EGOROV2} given \autoref{lem:SUFFICE}]
Given $T \in \mathbb{N}$ and two operators $U$ and $A$, we introduce the shorthand
\begin{equation}\label{SAMENORM}
\Ad^T(U)(A) = U^{T} A (U^*)^{T}
\end{equation}
for the $T$-fold conjugation of $A$ by $U$. To keep track of the remainders we
henceforth denote $R_N$ in the statement of \autoref{lem:SUFFICE} by $R_N^{(1)}$. Invoking the assumption \eqref{eqn:EXPGROWTH}  that  $\| \tilde{\chi}^T\|_{C^2}^2 = \ocal(e^{2|T|\delta_0})$, \autoref{lem:SUFFICE} reads
\begin{equation}\label{AdU}
\begin{dcases}
\Ad(U_{\chi,N}) \Pi_{h^N} M_F \Pi_{h^N} = \Pi_{h^N} M_{F \circ \tilde{\chi}} \Pi_{h^N} + R_N,\\
\frac{1}{d_N} \Tr[R_N^*R_N] = \ocal\bigg(\frac{1}{N}\| F\|_{C^2}^2 e^{2\delta_0}\bigg).
\end{dcases}
\end{equation}
We now iterate the conjugation.  
Conjugating a second time by $U_{\chi,N}$ yields two terms:
\begin{equation}\label{eqn:CONJTWICE1}
\Ad^2(U_{\chi,N}) \Pi_{h^N} M_F \Pi_{h^N} = \Ad(U_{\chi,N})\Pi_{h^N} M_{F \circ \tilde{\chi}} \Pi_{h^N} + \Ad(U_{\chi,N})R_N^{(1)}.
\end{equation}
It follows from \autoref{lem:SUFFICE} (with $M_F$ replaced by $M_{F \circ \tilde{\chi}}$) that the first  term on the right-hand side of \eqref{eqn:CONJTWICE1} equals
\begin{equation}\label{eqn:CONJTWICE2}
\begin{dcases}
 \Ad(U_{\chi,N})\Pi_{h^N} M_{F \circ \tilde{\chi}} \Pi_{h^N} = \Pi_{h^N} M_{F \circ \tilde{\chi}^2} \Pi_{h^N} + \tilde{R}^{(2)}_N, \\
\frac{1}{d_N}\Tr[(\tilde{R}_N^{(2)})^* \tilde{R}_N^{(2)}] = \ocal\bigg(\frac{1}{N}\| F \circ \tilde{\chi}\|_{C^2}^2 e^{2\delta_0}\bigg) = \ocal\bigg(\frac{1}{N}\| F\|_{C^2}^2 e^{4\delta_0}\bigg).
\end{dcases}
\end{equation}
In the error estimate we again made use of the exponential growth assumption \eqref{eqn:EXPGROWTH}.

The unitarity of $U_{\chi,N}$ implies that the second term $\Ad(U_{\chi,N}) R_N^{(1)}$ in \eqref{eqn:CONJTWICE2} satisfies
\begin{equation} \label{eqn:CONJTWICE3}
\Tr[(\Ad(U_{\chi,N})R_N^{(1)})^*\Ad(U_{\chi,N})R_N^{(1)}] = \Tr[(R_N^{(1)})^* R_N^{(1) }] =  \ocal\bigg(\frac{1}{N}\| F\|_{C^2}^2 e^{2\delta_0}\bigg).
\end{equation}
Combining \eqref{eqn:CONJTWICE1} and \eqref{eqn:CONJTWICE2} gives
\begin{equation}\label{eqn:CONJTWICE4}
\Ad^2(U_{\chi,N}) \Pi_{h^N} M_F \Pi_{h^N} = \Pi_{h^N} M_{F \circ \tilde{\chi}^2} \Pi_{h^N} + \tilde{R}^{(2)}_N + \Ad(U_{\chi,N})R_N^{(1)}.
\end{equation}
Set
\begin{equation}\label{eqn:CONJTWICE5}
R_N^{(2)} := \tilde{R}_N^{(2)} + \Ad(U_{\chi,N})R_N^{(1)},
\end{equation}
then \eqref{eqn:CONJTWICE2} and \eqref{eqn:CONJTWICE3} imply
\begin{align}
\frac{1}{d_N} \Tr[(R_N^{(2)})^* R_N^{(2)}] &\le \frac{2}{d_N} \Tr[(\tilde{R}_N^{(2)})^*\tilde{R}_N^{(2)} + (R_N^{(1)})^*R_N^{(1)}]\\
&= 2 \left(\ocal\bigg(\frac{1}{N}\| F\|_{C^2}^2 e^{4\delta_0}\bigg) + \ocal\bigg(\frac{1}{N}\| F\|_{C^2}^2 e^{2\delta_0}\bigg)\right)\\
&= 3 \ocal\bigg(\frac{1}{N}\| F\|_{C^2}^2 e^{4\delta_0}\bigg). \label{eqn:CONJTWICE6}
\end{align}
The statement of \autoref{prop:EGOROV2} with $T = 2$ is proved thanks to \eqref{eqn:CONJTWICE4}, \eqref{eqn:CONJTWICE5} and \eqref{eqn:CONJTWICE6}.

The calculation is similar when $\Ad(U_{\chi,N})$ is iterated $T$ times. By a similar stationary phase computation presented in the subsequent section, it is easy to see that on the $T$th
iterate, 
we pick up the leading order term:
\begin{equation}
\begin{dcases}
\Ad(U_{\chi,N})\Pi_{h^N} M_{F \circ \tilde{\chi}^{T-1}} \Pi_{h^N} = \Pi_{h^N} M_{F \circ \tilde{\chi}^{T}} \Pi_{h^N} + \tilde{R}^{(T)}_N,\\
\frac{1}{d_N}\Tr[(\tilde{R}_N^{(T)})^* \tilde{R}_N^{(T)}] = \ocal\bigg(\frac{1}{N}\| F\|_{C^2}^2 e^{2\delta_0|T|}\bigg).
\end{dcases}
\end{equation}
We also have to conjugate the $(T-1)$  `old' remainders from the $(T-1)$st iterate:
\begin{equation}
\Ad(U_{\chi,N}) \tilde{R}_N^{(T-1)} + \Ad^2(U_{\chi,N}) \tilde{R}_N^{(T-2)} + \Ad^3(U_{\chi,N}) \tilde{R}_N^{(T-3)} + \dotsb + \Ad^{T-1}(U_{\chi,N})\tilde{R}_N^{(1)}.
\end{equation}
The Hilbert-Schmidt norm of $\tilde{R}_N^{(\ell)}$ does not change under conjugation by $U_{\chi,N}$. Therefore, the combined remainder term
\begin{equation}
R_N^{(T)} := \tilde{R}^{(T)}_N + \Ad(U_{\chi,N}) \tilde{R}_N^{(T-1)} + \Ad^2(U_{\chi,N}) \tilde{R}_N^{(T-2)} +  \dotsb + \Ad^{T-1}(U_{\chi,N})\tilde{R}_N^{(1)}
\end{equation}
at the $T$th stage of the iterate has the estimate
\begin{equation}
\frac{1}{d_N} \Tr[(R_N^{(T)})^*R_N^{(T)}] \le \frac{T}{d_N} \sum_{\ell = 1}^T \Tr[(R_N^{(\ell)})^*R_N^{(\ell)}] = T \sum_{\ell = 1}^T \ocal\bigg(\frac{1}{N}\| F\|_{C^2}^2 e^{2\delta_0|\ell|}\bigg).
\end{equation}
Replacing each $e^{2\delta_0|\ell|}$ in the above sum by $e^{2\delta_0|T|}$ for $\ell = 1, 2, \dotsc, T$ completes the proof of \autoref{prop:EGOROV2} assuming \autoref{lem:SUFFICE}.
\end{proof}


\subsection{Proof of \autoref{lem:SUFFICE} via stationary phase computation}

Let
\begin{equation}
\tilde{L}_N := U_{\chi, N} \Pi_{h^N}   M_F \Pi_{h^N} U_{\chi, N}^* \quad \text{and} \quad L_N := \Pi_{h^N} M_{F \circ \tilde{\chi}} \Pi_{h^N},
\end{equation}

From the definition \eqref{eqn:QUANTIZE1} of Toeplitz quantization, the conjugated operator has the form
\begin{equation} \label{CONJ}
\tilde{L}_N = \Pi_{h^N} \sigma_N T_{\tilde{\chi}} \Pi_{h^N} M_F \Pi_{h^N}  T_{\tilde{\chi}^{-1}}\bar{\sigma}_N \Pi_{h^N}.
\end{equation}
Next, insert the identity operator $\operatorname{Id} = T_{\tilde{\chi}^{-1}} T_{\tilde{\chi}}$ between the operators $\Pi_{h^N}$ and $M_F$ in the above expression. Note that $T_{\tilde{\chi}} F T_{\tilde{\chi}^{-1}} =  F \circ \tilde{\chi}$. Hence, the expression becomes
\begin{equation}\label{eqn:GOAL}
\tilde{L}_N = \Pi_{h^N} \sigma_N  \Pi_{h^N}^{\tilde{\chi}} M_{F \circ \tilde{\chi}} \Pi_{h^N}^{\tilde{\chi}} \bar{\sigma}_N \Pi_{h^N}.
\end{equation}
where
$\Pi_{h^N}^{\tilde{\chi}} := T_{\tilde{\chi}} \Pi_{h^N} T_{\tilde{\chi}^{-1}}$ is the operator
with Schwartz kernel $\Pi_{h^N}^{\tilde{\chi}}(x,y) = \Pi_{h^N}(\chi(\tilde{x}), \chi(\tilde{y}))$.

In the notation  \eqref{eqn:SUFFICE}, 
\begin{equation}
R_N = \tilde{L}_N - L_N =  \Pi_{h^N} \big(\sigma_N  \Pi_{h^N}^{\tilde{\chi}} M_{F \circ \tilde{\chi}} \Pi_{h^N}^{\tilde{\chi}} \bar{\sigma}_N -  M_{F \circ \tilde{\chi}}\big) \Pi_{h^N}.
\end{equation}
Evidently,
\begin{equation} \label{CANCELS}
\Tr [R_N^* R_N] = \Tr [\tilde{L}_N^* \tilde{L}_N] - 2 \Tr [\tilde{L}_N L_N] + \Tr [L_N^* L_N].
\end{equation}
We evaluate each term asymptotically by stationary phase with remainder
and add the terms. \autoref{lem:SUFFICE} follows from:

\begin{lem}\label{lem:WEYL}
We have
\begin{equation} \label{TRLN2} \frac{1}{d_N}  \Tr [L_N^* L_N] = \int_M |F \circ  \tilde{\chi}|^2 \,d V + \ocal\bigg(\frac{1}{N}\| F\|_{C^2}^2 e^{2\delta_0}\bigg). \end{equation}
Moreover, 
\begin{equation}
\frac{1}{d_N} \Tr [\tilde{L}_N^*\tilde{L}_N] = \frac{1}{d_N} \Tr [\tilde{L}_N^*L_N] +  \ocal\bigg(\frac{1}{N}\| F\|_{C^2}^2 e^{2\delta_0}\bigg) = \frac{1}{d_N}\Tr[L_N^*L_N] + \ocal\bigg(\frac{1}{N}\| F\|_{C^2}^2 e^{2\delta_0}\bigg).
\end{equation}
In particular, thanks to \eqref{CANCELS} we have
\begin{equation}
\frac{1}{d_N} \Tr[R_N^*R_N] = \ocal \bigg(\frac{1}{N}\| F\|_{C^2}^2 e^{2\delta_0}\bigg).
\end{equation}
\end{lem}

The first statement \eqref{TRLN2} is the well-known \szego limit formula
with remainder. Since $\tilde{\chi}$ is symplectic it may be removed
from $F \circ \tilde{\chi}$ in the integral. The leading order term is calculated in \cite{BG81} using the homogeneous calculus of Toeplitz operators. The semi-classical calculation
and the remainder estimate may be calculated by the method below.

For the rest of the Appendix, we calculate the most difficult of the three terms, namely $d_N^{-1}\Tr [\tilde{L}_N^*\tilde{L}_N]$, asymptotically to leading order by the method
of stationary phase for oscillatory integrals with complex phases of positive
type (\cite{H90}, Theorem
7.7.5). We use the remainder estimate from that theorem. The calculations
of the other two terms are similar and therefore omitted.

  All three traces in \eqref{CANCELS} have the same leading order term
  \eqref{TRLN2},  and so the leading term cancels when taking
the sum \eqref{CANCELS}. The cancellation between the `symbols' $\sigma_N$  and the Hessian determinants in the calculation of the leading order terms  \eqref{TRLN2} is guaranteed by unitarity of 
$U_{\chi, N}$ (see also \cite{Z97} for explicit calculation of the symbol).

From \eqref{eqn:GOAL}, we have
\begin{equation}\label{eqn:TRACE}
\frac{1}{d_N} \Tr[\tilde{L}_N^*\tilde{L}_N] = \frac{1}{d_N} \Tr \big[\Pi_{h^N} \bar{\sigma}_N  \Pi_{h^N}^{\tilde{\chi}} M_{\overline{F \circ \tilde{\chi}}} \Pi_{h^N}^{\tilde{\chi}} \sigma_N \Pi_{h^N}\sigma_N  \Pi_{h^N}^{\tilde{\chi}} M_{F \circ \tilde{\chi}} \Pi_{h^N}^{\tilde{\chi}} \bar{\sigma}_N \big].
\end{equation}
Note that we may drop the factor of $\Pi_{h^N}$ at the end when computing the trace. We use the shorthand
\begin{equation}
\quad \tilde{y}_j := \tilde{\chi}(y_j), \quad y_j \in X.
\end{equation}
Recall that $\sigma_N$ denotes multiplication by the symbol $\sigma_N$, and the \szego projectors have Schwartz kernels
\begin{align}
\Pi_{h^N}^{\tilde{\chi}}(y_1,y_2) &= \Pi_{h^N}(\tilde{y}_1,\tilde{y}_2),\\
\Pi_{h^N}(y_1,y_2) &=N\int_0^\infty\!\int_{S^1} e^{iN[-\theta + t \psi(r_\theta y_1, y_2)]} s(r_\theta y_1,y_2,Nt)\,d\theta dt. 
\end{align}
The last equality is the \BSj parametrix introduced in \autoref{sec:BSj}. Using Schwartz kernels, the trace \eqref{eqn:TRACE} can be written as the following oscillatory integral
\begin{align}
\frac{1}{d_N} \Tr[\tilde{L}_N^*\tilde{L}_N] &= \frac{1}{d_N}\int_X (\tilde{L}_N^*\tilde{L}_N)(x,x)\,dx\\
&= \frac{1}{d_N}\int_X\left(N^6 \int_{X^5 \times (S^1)^6 \times (\R_+)^6} A(x,\vec{y},\vec{\theta},\vec{t}) e^{iN \Psi(x,\vec{y},\vec{\theta},\vec{t})}\,d\vec{t}d\vec{\theta}d\vec{y}\right)dx,
\end{align}
where
\begin{equation}
\vec{y} = (y_1, \dotsc, y_5) \in X^5, \quad \vec{\theta} = (\theta_1, \dotsc, \theta_6) \in (S^1)^6, \quad \vec{t} = (t_1, \dotsc, t_6) \in (\R_+)^6
\end{equation}
and the amplitude and phase function are given by
\begin{align}
A &= s(r_{\theta_1} x,y_1,t_1N) \bar{\sigma}_N(y_1) s(r_{\theta_2}\tilde{y}_1,\tilde{y}_2,t_2N) \overline{F(\tilde{y}_2)}s(r_{\theta_3}\tilde{y}_2,\tilde{y}_3,t_3N)\sigma_N(y_3)\\
&\quad \times s(r_{\theta_4}y_3,y_4,t_4N) \sigma_N(y_4) s(r_{\theta_5}\tilde{y}_4,\tilde{y}_5,t_5N) F(\tilde{y}_5) s(r_{\theta_6}\tilde{y}_5, \tilde{x},t_6N) \bar{\sigma}_N(x),\\
\Psi &= t_1\psi(r_{\theta_1}x,y_1) - \theta_1 + t_2\psi(r_{\theta_2}\tilde{y}_1, \tilde{y}_2) - \theta_2 + t_3\psi(r_{\theta_3}\tilde{y}_2,\tilde{y}_3) - \theta_3\\
&\quad  + t_4\psi(r_{\theta_4}y_3,y_4) - \theta_4 + t_5 \psi(r_{\theta_5}\tilde{y}_4,\tilde{y}_5) - \theta_5 + t_6 \psi(r_{\theta_6}\tilde{y}_5,\tilde{x}).
\end{align}
The functions $s$ and $\psi$ come from the \BSj parametrix \eqref{eqn:PARAMETRIX}, and $\sigma_N$ comes from the quantization formula \eqref{eqn:QUANTIZE1}.

The method of stationary phase is used to compute the inner integral. The off-diagonal exponential decay estimate \eqref{offdiag} for the Bergman kernel allows us to localize the $X^5$-space integral to the region $\{d(y_j,y_k) < N^{-1/3}\}$ and absorb the error in the remainder estimate for $R_N$. To locate the critical points of the phase function $\Psi$, recall from \eqref{eqn:SMALL PHASE} that the function $\psi$ has the form
\begin{equation}
\psi(x,y) = \frac{1}{i}\bigg(1  - \Lambda(x,y)\bigg) \quad \text{with} \quad \Lambda(x,y) := e^{-\frac{\phi(z_1)}{2} - \frac{\phi(z_2)}{2} + \phi(z_1,\bar{z}_2)}e^{i(\tau_1 - \tau_2)},
\end{equation}
from which it follows
\begin{equation}
\psi(r_\theta x,y) = \frac{1}{i} \bigg( 1- e^{i\theta}\Lambda(x,y)\bigg).
\end{equation}
Therefore,
\begin{equation}
D_{t_1} \Psi = \psi(r_{\theta_1}x,y_1) = 0 \iff 1 = e^{i\theta_1}\Lambda(x,y_1).
\end{equation}
The Schwarz inequality shows that a real critical point exists if and only if $x = y_1$. Similar computations for $D_{t_j}\Psi$ demand that $\tilde{y}_1 = \tilde{y}_2 = \tilde{y}_3$, $y_3 = y_4$, and $\tilde{y}_4 = \tilde{y}_5 = \tilde{x}$. The real critical point of $\Psi$ must therefore satisfy
\begin{equation}\label{eqn:CRIT1}
x = y_1 = y_2 = y_3 = y_4 = y_5.
\end{equation}

Consider now the $\theta_1$ derivative:
\begin{equation}
D_{\theta_1} \Psi = - t_1e^{i\theta_1}\Lambda (x,y_1) - 1 = 0 \iff 1 = -t_1e^{i\theta_1}\Lambda(x,y_1).
\end{equation}
From the constraint \eqref{eqn:CRIT1}, we must have $x = (z_1,\tau_1) = (z_2,\tau_2) = y_1$, so $\Lambda(x,y_1) = 1$. It follows that $t_1 - -1$ and $\theta_1 = 0$. Similar computations for $D_{\theta_j}\Psi$ show that the real critical point of $\Psi$ satisfies
\begin{equation}\label{eqn:CRIT2}
\theta_1 = \dotsb = \theta_6 = 0 \quad \text{and} \quad t_1 = \dotsb = t_6  = -1.
\end{equation}

Finally, we claim that $D_{y_j}\Psi$ automatically vanishes at the points satisfying \eqref{eqn:CRIT1} and \eqref{eqn:CRIT2}. Indeed, at the critical point we have
\begin{equation}
D_{y_1} \Psi\bigg|_{\substack{x = y_1 = \dotsb = y_5 \\ \theta_j = 0 \\ t_j = -1}} = -D_{y_1} \psi(x,y_1)|_{y_1 = x} - D_{y_1}\psi(\tilde{y}_1,\tilde{y}_2)|_{y_2 = y_1 = x}.
\end{equation}
Recall, however, that along the diagonal of $X \times X$ we have
\begin{equation}
d_1 \psi = -d_2\psi = \frac{1}{i} d\rho|_X = \alpha,
\end{equation}
where $\alpha$ is the contact form. Here $d_j$ refers to the derivative with respect to the $j$th slot of $\psi(\cdot, \cdot)$. The assumption that $\chi$ lifts to a contact transformation, that is, $\tilde{\chi}^*\alpha = \alpha$, implies
\begin{equation}
-D_{y_1} \psi(x,y_1)|_{y_1 = x} - D_{y_1}\psi(\tilde{y}_1,\tilde{y}_2)|_{y_2 = y_1 = x} = \alpha(x) - \frac{1}{i}d\rho(\tilde{\chi}(x)) = \alpha(x) - \tilde{\chi}^*\bigg(\frac{1}{i}d\rho\bigg)(x) = 0.
\end{equation}
Similar computations for $D_{y_j} \Psi$ show that the real critical points of $\Psi$ are completely given by \eqref{eqn:CRIT1} and \eqref{eqn:CRIT2}.

It is straightforward to verify that the Hessian at the critical point is a block matrix of the form
\begin{equation}
\Hess \Psi(x) = \begin{bmatrix} 
D_{\vec{t}\vec{t}}\Psi = 0 & D_{\vec{t}\vec{\theta}} \Psi = -\mathrm{Id} & D_{\vec{t}\vec{1}}\Psi & D_{\vec{t}\vec{2}}\Psi\\
D_{\vec{\theta}\vec{t}} \Psi = -\mathrm{Id} & D_{\vec{\theta}\vec{\theta}} \Psi = i \cdot \mathrm{Id} & D_{\vec{\theta}\vec{1}}\Psi & D_{\vec{\theta}\vec{2}}\Psi\\
D_{\vec{1}\vec{t}}\Psi & D_{\vec{1}\vec{\theta}}\Psi & D_{\vec{1}\vec{1}} \Psi & D_{\vec{1}\vec{2}}\Psi \\
D_{\vec{2}\vec{t}}\Psi & D_{\vec{2}\vec{\theta}}\Psi & D_{\vec{2}\vec{1}} \Psi & D_{\vec{2}\vec{2}}\Psi
\end{bmatrix}
\end{equation}
with
\begin{equation}
D_{\vec{t}\vec{1}}\Psi = - D_{\vec{t}\vec{2}}\Psi = \begin{bmatrix}
\alpha(x) & 0 & 0 & 0 & 0\\
-\alpha(x) & \alpha(x) &0 & 0 & 0\\
0 & \ddots & \ddots & 0 & 0\\
0 & 0 & \ddots & \ddots & 0 \\
0 & 0 & 0 & -\alpha(x) & \alpha(x)\\
0 & 0 & 0 & 0 & -\alpha(x) \end{bmatrix}
= - (D_{\vec{2}\vec{t}}\Psi)^t = (D_{\vec{1}\vec{t}}\Psi)^t,
\end{equation}
\begin{equation}
D_{\vec{\theta}\vec{1}}\Psi = - D_{\vec{\theta}\vec{2}}\Psi = \begin{bmatrix}
-i\alpha(x) & 0 & 0 & 0 & 0\\
i\alpha(x) & -i\alpha(x) & 0 & 0 & 0 \\
0 & \ddots & \ddots &0 & 0 \\
0 & 0 & \ddots & \ddots & 0 \\ 
0 & 0 & 0 & i\alpha(x) & -i\alpha(x)\\
0 & 0 & 0 & 0 & i\alpha(x)
\end{bmatrix} = - (D_{\vec{2}\vec{\theta}})^t = (D_{\vec{2}\vec{\theta}}\Psi)^t,
\end{equation}
\begin{equation}
D_{\vec{1}\vec{1}}\Psi = \begin{bmatrix} -d\alpha(x) &d\alpha(x)&0 &0 & 0\\
d\alpha(x)& \ddots & \ddots & 0 & 0\\
0 & \ddots & \ddots & \ddots & 0 \\
0 & 0 & \ddots & \ddots & d\alpha(x) \\
0& 0 & 0 &d\alpha(x)& -d\alpha(x)  \end{bmatrix} = D_{\vec{2}\vec{2}}\Psi.
\end{equation}
This Hessian matrix is invertible by the Schur complement formula (recall that $-id\rho = \alpha$ is non-vanishing in a neighborhood of $X$). The method of stationary phase shows that the Schwartz kernel $(\tilde{L}_N^*\tilde{L}_N)(x,x)$ along the diagonal has the expansion
\begin{multline}\label{STNPHSEXP}
(\tilde{L}_N^*\tilde{L}_N)(x,x) \sim \frac{N^6}{(N^{12 + 10m}\det \Hess\Psi(x))^{1/2}} \sum_{j,k,\ell,p,q,u,v \ge 0}N^{6m-j-k-\ell-p-q -u -v} \\
\times L_j\Big(s_k(x,x)s_\ell(x,x)s_p(\tilde{x},\tilde{x})s_q(\tilde{x},\tilde{x}) s_u(\tilde{x},\tilde{x})s_v(\tilde{x},\tilde{x})\lvert \sigma_N(x)\rvert^4 \lvert F(\tilde{x})\rvert^2\Big),
\end{multline}
where $L_j$ are differential operators of order at most $2j$ that can be explicitly expressed in terms of $s_k$ and the Hessian (see~\cite{H90}).

Observe that the leading order term (obtained from the above expression by setting $j = k = \dotsb = v = 0$) is of order $N^6(N^{12 + 10m})^{-1/2} N^{6m} = N^m$. The symbol  $\sigma_N$ is constructed to make $U_{\chi,N}$ unitary, i.e., $U_{\chi,N}^* U_{\chi,N} = \Pi_{h^N}$, and by taking the symbol
of this equation it follows that
\begin{equation}\label{LEADCANCEL}
(\det \Hess\Psi(x))^{-1/2} s_0(x,x)^2  s_0(\tilde{x},\tilde{x})^4|\sigma_N(x)|^4 = 1.
\end{equation}
Indeed, if we set $F \equiv 1$ so that $M_F = \operatorname{Id}$, then $\tilde{L}_N^*\tilde{L}_N = U_{\chi,N}U_{\chi,N}^*U_{\chi,N}U_{\chi,N}^* = \operatorname{Id}$. The identity \eqref{LEADCANCEL} follows from plugging this particular choice of $F$ into \eqref{STNPHSEXP}. Therefore, after dividing by $d_N \sim N^{m}$ (for $N$ large enough), the leading order term of $d_N^{-1}\Tr[\tilde{L}_N^*\tilde{L}_N]$ is of order $0$, and is equal to $\int |F (\tilde{x})|^2 = \int |F \circ \tilde{\chi}|^2$, which agrees with \eqref{TRLN2}. The second order term (cf.\ \cite[Theorem 7.7.5]{H90}) of $\tilde{L}_N^*L_N(x,x)$ is bounded above in sup norm by
\begin{align}
C \sum_{|\alpha| \le 2} \Big\| D^\alpha \Big((\det \Hess\Psi(x))^{-\frac{1}{2}}s_0(x,x)^2s_0(\tilde{\chi}(x),\tilde{\chi}(x))^4\lvert \sigma_N(x) \rvert^4 &\lvert F \circ \tilde{\chi}(x) \rvert^2\Big)\Big\|_\infty \\
&= C \sum_{|\alpha| \le 2} \Big\| D^\alpha \lvert F \circ \tilde{\chi}(x) \rvert^2 \Big\|_\infty\\
&\le C \bigg(\sum_{|\alpha| \le 2} \Big\| D^\alpha \lvert F \circ \tilde{\chi}(x) \rvert \Big\|_\infty\bigg)^2\\
&\le C \| F\|_{C^2}^2 e^{2\delta_0}.
\end{align}
for some constant $C$ that depends on a fixed number of derivatives of the phase function $\Psi$ (and hence on $\tilde{\chi}$) but  is otherwise independent of $N$. Dividing through by $d_N \sim N^m$ yields the desired error estimate $\ocal(N^{-1} \| F\|_{C^2}^2 e^{2\delta_0})$. This completes the computation for $\tilde{L}_N$.


\begin{thebibliography}{HHHH}


\bibitem[BP]{BP13} L. Barreira\ and\ Y. Pesin, \textit{Introduction to smooth ergodic theory, Graduate Studies in Mathematics}, 148, American Mathematical Society, Providence, RI (2013). MR3076414




 

\bibitem[BG]{BG81} L. Boutet de Monvel\ and\ V. Guillemin, \emph{The spectral theory of Toeplitz operators, Annals of Mathematics Studies}, 99, Princeton University Press, Princeton, NJ (1981). MR0620794


\bibitem[BS]{BS75} L. Boutet de Monvel\ and\ J. Sj\" ostrand, Sur la singularit\'e des noyaux de Bergman et de \szego, Journ\'ees: \'Equations aux D\'eriv\'ees Partielles de Rennes (1975), 123--164. Ast\'erisque, 34-35, Soc. Math. France, Paris. MR0590106.

\bibitem[BouR]{BouR02}  A. Bouzouina and D. Robert,
Uniform semiclassical estimates for the propagation of quantum observables. 
Duke Math. J. 111 (2002), no. 2, 223-252. 

\bibitem[BW]{BW10}  K. Burns\ and\ A. Wilkinson, On the ergodicity of partially hyperbolic systems, Ann. of Math. (2) {\bf 171} (2010), no.~1, 451--489. MR2630044


\bibitem[De]{Delin98} H. Delin, Pointwise estimates for the weighted Bergman projection kernel in $\C^n$, using a weighted $L^2$ estimate for the $\overline\partial$ equation, Ann. Inst. Fourier (Grenoble)  48 (1998), no.~4, 967--997. MR1656004



\bibitem[DeV]{DeV11} K. Dekimpe and K. Verheyen, 
Constructing infra-nilmanifolds admitting an Anosov diffeomorphism. 
Adv. Math. 228 (2011), no. 6, 3300-3319. 


\bibitem[DW]{DW03} D. Dolgopyat\ and\ A. Wilkinson, Stable accessibility is $C^1$ dense, Ast\'erisque No. 287 (2003), xvii, 33--60. MR2039999

\bibitem[FT]{FT15}  F. Faure and M.  Tsujii, 
pre-quantum transfer operator for symplectic Anosov diffeomorphism.
Ast\'erisque No. 375 (2015).



\bibitem[GH]{GriffithsHarris78} P. Griffiths\ and\ J. Harris, {\it Principles of algebraic geometry}, Wiley-Interscience, New York, 1978. MR0507725


\bibitem[Ha]{H15} X. Han, Small scale quantum ergodicity in negatively curved manifolds, Nonlinearity {\bf 28} (2015), no.~9, 3263--3288. MR3403398.


\bibitem[HR]{HR16} H. Hezari\ and\ G. Rivi\`ere, $L^p$ norms, nodal sets, and quantum ergodicity, Adv. Math. {\bf 290} (2016), 938--966. MR3451943

\bibitem[H\"o]{H90} L. H\"ormander, \emph{The analysis of linear partial differential operators. I}, reprint of the second (1990) edition [Springer, Berlin; MR1065993 (91m:35001a)], Classics in Mathematics, Springer, Berlin, 2003. MR1996773


\bibitem[Ka]{Ka79} A. Katok, Bernoulli diffeomorphisms on surfaces, Ann. of Math. (2) {\bf 110} (1979)no.~3. MR0554383

\bibitem[Ke]{Kel10} D.  Kelmer, Arithmetic quantum unique ergodicity for symplectic linear maps of the multidimensional torus, Ann. of Math. (2) 171 (2010), no. 2, 815--879. MR2630057

 

\bibitem[L]{Lindholm01} N. Lindholm, Sampling in weighted $L^p$ spaces of entire functions in ${\C}^n$ and estimates of the Bergman kernel, J. Funct. Anal. {\bf 182} (2001), no.~2, 390--426. MR1828799

\bibitem[LMR]{LMR15} S. Lester, K. Matom\"{a}ki and M. Radziwi\l\l, Small scale distribution of zeros and mass of modular forms. (English summary) J. Eur. Math. Soc. (JEMS) 20 (2018), no. 7, 1595–1627. 

\bibitem[LR]{LR} S. Lester\ and\ Z. Rudnick, Small scale equidistribution of eigenfunctions on the torus, Comm. Math. Phys. {\bf 350} (2017), no.~1, 279--300. MR3606476

\bibitem[LuSh]{LuSh15}  Z. Lu and B.  Shiffman,  Asymptotic expansion of the off-diagonal Bergman kernel on compact \kahler manifolds. J. Geom. Anal. 25 (2015), no. 2, 761-782. 

\bibitem[MM]{MaMa12} X. Ma and G.  Marinescu, 
Berezin-Toeplitz quantization on Kaehler manifolds.
J. Reine Angew. Math. 662 (2012), 1-56. 

 
 \bibitem[Ma]{M16} K. Marin, $C^r$-density of (non-uniform) hyperbolicity in partially hyperbolic symplectic diffeomorphisms, Comment. Math. Helv. {\bf 91} (2016), no.~2, 357--396. MR3493375
 
 \bibitem[NV]{NV98}  S. Nonnenmacher\ and\ A. Voros, Chaotic eigenfunctions in phase space, J. Statist. Phys. {\bf 92} (1998), no.~3-4, 431--518. MR1649013






\bibitem[R]{R05} Z. Rudnick, On the asymptotic distribution of zeros of modular forms, Int. Math. Res. Not. {\bf 2005}, no.~34, 2059--2074. MR2181743


\bibitem[Sc1]{Sc06} R. Schubert, Upper bounds on the rate of quantum ergodicity, Ann. Henri Poincar\'e {\bf 7} (2006), no.~6, 1085--1098. MR2267060

\bibitem[Sc2]{Sc08} R. Schubert, On the rate of quantum ergodicity for quantised maps, Ann. Henri Poincar\'e {\bf 9} (2008), no.~8, 1455--1477. MR2465731


\bibitem[ShZ1]{ShZ99} B. Shiffman\ and\ S. Zelditch, Distribution of zeros of random and quantum chaotic sections of positive line bundles, Comm. Math. Phys. {\bf 200} (1999), no.~3, 661--683. MR1675133

\bibitem[ShZ2]{ShZ02} B. Shiffman and S. Zelditch, Asymptotics of almost holomorphic sections of ample line bundles on symplectic manifolds, J. Reine Angew. Math., 544 (2002), pp. 181-222. MR1887895








\bibitem[Z1]{Z97} S. Zelditch, Index and dynamics of quantized contact transformations, Ann. Inst. Fourier (Grenoble) {\bf 47} (1997), no.~1, 305--363. MR1437187


\bibitem[Z2]{Z98} S. Zelditch, \szego kernels and a theorem of Tian, Internat. Math. Res. Notices {\bf 1998}, no.~6, 317--331. MR1616718



\end{thebibliography}
\end{document}